\documentclass[aop]{imsart}

\RequirePackage{amsthm,amsmath,amsfonts,amssymb}
\RequirePackage[numbers]{natbib}
\usepackage{mathrsfs}

\startlocaldefs
\newtheorem{theorem}{Theorem}
\newtheorem{proposition}[theorem]{Proposition}
\newtheorem{definition}[theorem]{Definition}
\newtheorem{lemma}[theorem]{Lemma}

\def\build#1_#2^#3{\mathrel{
\mathop{\kern 0pt#1}\limits_{#2}^{#3}}}

\def\w{\mathrm{w}}

\def\t{\mathcal{T}}

\def\W{\mathcal{W}}
\def\S{\mathcal{S}}

\def\N{\mathbb{N}}

\def\D{\mathbb{D}}
\def\P{\mathbb{P}}

\def\E{\mathbb{E}}

\def\R{\mathbb{R}}
\def\z{\mathcal{Z}}
\def\y{\mathcal{Y}}
\def\n{\mathcal{N}}

\def\ve{{\varepsilon}}
\def\la{\longrightarrow}

\def\ov{\overline}

\def\dd{\mathrm{d}}
\def\wh{\widehat}
\def\wt{\widetilde}
\def\tr{\mathrm{tr}}

\def\v{\mathcal{V}}

\def\XX{\mathcal{X}}

\def\bn{\mathbf{n}}

\def\bm{\mathbf{m}}
\def\bx{\mathbf{x}}

\def\rems{\noindent{\bf Remarks. }}

\endlocaldefs

\begin{document}

\begin{frontmatter}
\title{The Markov property of local times \\
of Brownian motion indexed by the Brownian tree}
\runtitle{Markov property of local times}

\begin{aug}
\author[A]{\fnms{Jean-Fran\c cois}~\snm{Le Gall}\ead[label=e1]{jean-francois.le-gall@universite-paris-saclay.fr}},
\address[A]{D\'epartement de Math\'ematiques d'Orsay, Universit\'e Paris-Saclay\printead[presep={,\ }]{e1}}

\end{aug}

\begin{abstract}
We consider the model of Brownian
motion indexed by the Brownian tree, which has appeared in 
a variety of different contexts in probability, statistical physics and combinatorics.
For this model, the total occupation measure is known to have 
a continuously differentiable density $(\ell^x)_{x\in\R}$ and we write 
$(\dot\ell^x)_{x\in\R}$ for its derivative. Although the process $(\ell^x)_{x\geq 0}$
is not Markov, we prove that the pair $(\ell^x,\dot\ell^x)_{x\geq 0}$
is a time-homogeneous Markov process. We also establish a
similar result for the local times of one-dimensional super-Brownian 
motion. Our methods rely on the excursion theory for Brownian
motion indexed by the Brownian tree.

\end{abstract}



\end{frontmatter}


\section{Introduction}

The Ray-Knight theorems, which give the Markov property 
of the process
of local times of linear Brownian motion in the space variable, at certain particular stopping
times, are some of the most famous and useful results about
Brownian motion. The goal of the present work is to discuss a similar Markov property
of local times for
the model of branching Brownian motion which we call Brownian
motion indexed by the Brownian tree. Here the Brownian tree 
is conveniently described as the random continuous tree $\t$ coded 
by a Brownian excursion  under the It\^o measure, and may
also be viewed as a free version of Aldous' Continuum Random Tree
(the word ``free'' means that the volume of the tree is not fixed). The tree
$\t$ is equipped with a volume measure $\mathrm{Vol}$ whose total
mass is the duration $\sigma$ of the excursion coding $\t$. 
Given the Brownian tree $\t$, we can consider Brownian motion
indexed by $\t$, which we denote by $(V_a)_{a\in\t}$. We view
$V_a$ as a label assigned to the  ``vertex'' $a$ of the tree, in such a way
that the label of the root is $0$, and labels evolve like linear
Brownian motion along line segments of $\t$. The total occupation
measure of Brownian motion indexed by $\t$ is the measure $\mathcal{Y}$ defined by
$$\langle \mathcal{Y},g\rangle = \int g(V_a)\,\mathrm{Vol}(\dd a),$$
for every nonnegative measurable function $g$ on $\R$. We write
$\N_0$ for the ($\sigma$-finite) measure under which $\t$
and $(V_a)_{a\in\t}$ are defined.

Let us emphasize that the pair $(\t,(V_a)_{a\in\t})$ plays an important role 
in a number of different areas of probability theory, combinatorics or
statistical physics. In particular, this pair is the key ingredient of the Brownian snake
construction of super-Brownian motion \cite{Zurich}. When conditioned on having
a total volume equal to $1$ (this just means that the coding Brownian excursion
is normalized), $\t$ becomes Aldous' Continuum Random Tree also known as the CRT \cite{Ald}, up to
an unimportant scaling factor $2$, and $\mathcal{Y}$ then corresponds
to the random measure called ISE \cite{Al2}. Both the CRT and ISE appear in
a number of combinatorial asymptotics for discrete trees
possibly equipped with labels (see e.g.~\cite{BM0,BMJ,probasur}). Other applications, involving 
multidimensional versions of $(V_a)_{a\in\t}$, include
interacting particle systems (see e.g.~\cite{BCL}) and models of statistical physics \cite{DS,HaS}. 
More recently, the pair $(\t,V_a)_{a\in\t})$ has been used as the
basic building block for the construction of the models 
of random geometry that arise as scaling limits of large random planar
maps (see in particular \cite{Uniqueness,Mie}).

The measure
$\mathcal{Y}$ has ($\N_0$ a.e.) a continuous density denoted by $(\ell^x,x\in \R)$,
and we call $\ell^x$ the local time of $(V_a)_{a\in\t}$ at level 
$x$. The function $x\mapsto \ell^x$ is even continuously differentiable on $\R$. The
latter property is proved
in the recent paper \cite{CM}, and a slightly weaker statement 
had been obtained earlier in \cite{BMJ} (both \cite{BMJ} and \cite{CM}
deal with ISE, but
a straightforward scaling argument then shows that the preceding
properties also hold for $\mathcal{Y}$ under $\N_0$).
As a matter of fact, the existence of a continuously differentiable density
for  $\mathcal{Y}$ under $\N_0$ can also be derived from older
results of Sugitani \cite{Sug}, which were concerned with 
(one-dimensional) super-Brownian motion (see Section \ref{sec:super} below). 
We write $\dot\ell^x$ for the derivative of the function $x\mapsto\ell^x$.

By analogy with the classical Ray-Knight theorems, one may ask
about the Markovian properties of the process $(\ell^x,x\in\R)$,
or simply of $(\ell^x,x\geq 0)$. Informally, it seems clear
that this process cannot be Markovian: To predict the future 
after time $x\geq 0$, one should at least know the value 
of the derivative $\dot\ell^x$, and not only $\ell^x$. The main result 
of the present work shows that the additional information
provided by the derivative suffices to get a Markov process.

\begin{theorem}
\label{mainT}
The process $((\ell^x,\dot\ell^x),x\geq 0)$ is time-homogeneous 
Markov under $\N_0$. Moreover the two processes $((\ell^x,\dot\ell^x),x\geq 0)$
and $((\ell^{-x},\dot\ell^{-x}),x\geq 0)$ are independent conditionally on $(\ell^0,\dot\ell^0)$.
\end{theorem}

A simple symmetry argument shows that $((\ell^x,\dot\ell^x),x\geq 0)$ 
and $((\ell^{-x},-\dot\ell^{-x}),x\geq 0)$ have the same distribution (in particular, the
law of $\dot\ell^0$ is symmetric). 
One may be puzzled by the fact that $\N_0$ is an infinite measure.
However, for every $\ve>0$, the event where $\ell^0\geq \ve$
has finite $\N_0$-measure (the distribution of $\ell^0$
under $\N_0$ has a density proportional to $\ell^{-5/3}$, cf. \cite[Corollary 3.1]{LGR2}) and the statement of the theorem 
can be formulated 
as well under the probability measure $\N_0(\cdot\mid \ell^0\geq \ve)$.

Let us now discuss an analog of Theorem \ref{mainT} for super-Brownian motion. We consider a 
one-dimensional super-Brownian motion $(\mathbf{X}_t)_{t\geq 0}$ 
with branching mechanism $\phi(u)=2u^2$ and initial value $\mathbf{X}_0=\alpha\delta_0$,
where $\alpha>0$ is a constant. Note that the choice of the multiplicative constant $2$
in the formula for $\phi$ is only for convenience and by scaling one could as well 
deal with the case $\phi(u)=c\,u^2$ for $c>0$. By results of Sugitani \cite[Theorem 4]{Sug}, the total occupation measure
$$\mathbf{Y}:=\int_0^\infty \mathbf{X}_t \,\dd t.$$
has (a.s.) a 
continuous density $(L^x)_{x\in\R}$ with respect to Lebesgue measure, and 
this density is continuously differentiable on $(0,\infty)$
and on $(-\infty,0)$. Let $\dot L^y$ stand for the derivative 
of the mapping $x\mapsto L^x$ at $y\not =0$. When $y=0$, both the right derivative
$\dot L^{0+}$ and the left derivative $\dot L^{0-}$ exist, and $\dot L^{0+}-\dot L^{0-}=-2\alpha$.
By convention, we set $\dot L^0=\dot L^{0+}$.

\begin{theorem}
\label{main-super}
The process $((L^x,\dot L^x),x\geq 0)$ is time-homogeneous 
Markov with the same transition kernel as the process $((\ell^x,\dot\ell^x),x\geq 0)$ of Theorem \ref{mainT}. Moreover the two processes $((L^x,\dot L^x),x\geq 0)$
and $((L^{-x},\dot L^{-x}),x\geq 0)$ are independent conditionally on $(L^0,\dot L^0)$.
\end{theorem}

By symmetry, the two processes $((L^x,\dot L^x),x\geq 0)$ and $((L^{-x},-\dot L^{(-x)-}),x\geq 0)$ have the same law,
where $L^{(-x)-}=L^{-x}$ except when $x=0$. In particular, the law of $\dot L^0+\alpha$ is symmetric. 
Theorem \ref{main-super} is derived by adapting the method of proof of Theorem \ref{mainT},
using the fact that the process $(\mathbf{X}_t)_{t\geq 0}$ can be constructed from a Poisson point measure 
with intensity $\alpha\N_0$ (see \cite[Chapter IV]{Zurich}). 

Let us explain the main ideas of the proof of Theorem \ref{mainT}. 
It is well known that the classical Ray-Knight theorems can be proved
by excursion theory, using in particular the independence of
excursions above and below a given level. Our proof of Theorem \ref{mainT} follows a similar approach,
but we now rely on the excursion theory developed in the article
\cite{ALG} for Brownian motion indexed by the Brownian tree.
Let us fix $h>0$ for definiteness. As in the classical setting, one is interested
in describing the connected components of
the set $\{a\in\t:V_a\not =h\}$ together with the distribution of the Brownian labels $V_a$ assigned 
to each connected component. Leaving aside the connected component
containing the root of $\t$, which is called the root component and plays a particular role, we call any such component (equipped with its labels) an excursion above or below $h$,
depending on the fact that labels are greater or smaller than $h$.
For any excursion above or below $h$, one can make sense of a quantity called the boundary size of the excursion, which measures
how many points of the closure of the component have a label equal to $h$. Moreover, 
the boundary sizes of the components other than the root component are exactly the jumps of a continuous-state branching process
with stable branching mechanism $\psi(\lambda)=\sqrt{8/3}\,\lambda^{3/2}$, which we denote by $(\XX^h_r)_{r\geq 0}$ and whose
initial value $\z_h$ is the so-called ``exit measure'' from $(-\infty,h)$ 
($\z_h$ corresponds to the boundary size of the root component). Roughly speaking, 
the results of \cite{ALG} imply that the excursions above or below $h$ are independent 
(and are also independent of the root component) conditionally on
their boundary sizes. The point in deriving Theorem \ref{mainT} is now to understand the conditional
distribution of the boundary sizes of excursions above level $h$ given the excursions 
below $h$ (and the root component). To this end, we first observe that the classical Lamperti representation allows us to
write $(\XX^h_r)_{r\geq 0}$ as a time change of a stable L\'evy process $U$ with no negative jumps
started at $U_0=\z_h$ and stopped 
at the time $T_0$ when it first hits $0$.
The boundary sizes of excursions also correspond to the jumps of this (stopped) L\'evy process. 
Distinguishing excursions above and below level $h$ amounts to assigning a label $+1$ or $-1$
to each of these jumps.  One can construct two independent L\'evy processes $U'$ and $U''$,
such that, on one hand, $U'_0=\z_h$ and the jumps of $U'$ are the jumps of $U$ with label $-1$,
on the other hand, $U''_0=0$ and the jumps of $U''$ are the jumps of $U$ with label $+1$
(in such a way that $U=U'+U''$, and the L\'evy measure of $U'$, or of $U''$, is half the L\'evy measure of $U$).
Finally, one can prove that the local time $\ell^h$ is equal to $T_0$ and moreover its derivative
$\dot\ell^h$ is equal to $2U''_{T_0}=-2U'_{T_0}$. From these observations and some
additional work, one gets that the conditional distribution of the
boundary sizes of excursions above $h$, knowing the excursions below $h$ and the root component,
is the distribution of jumps of the L\'evy process $U''$ conditioned to be equal to $\frac{1}{2}\dot\ell^h$
at time $\ell^h$, and this conditional distribution only depends on the pair $(\ell^h,\dot\ell^h)$. 
This leads to the desired Markov property. 

It is interesting to compare Theorem \ref{mainT} with the main result of \cite{LGR}, which gives
the distribution under $\N_0$ of the random process $(\mathscr{X}_x)_{x\geq 0}$ whose value at time $x\geq 0$
is the sequence of boundary sizes of connected components of $\{a\in\t:V_a>x\}$ 
in nonincreasing order
(these are the boundary sizes of excursions above level $x$, in the language of
the preceding paragraph). The process $(\mathscr{X}_x)_{x\geq 0}$ is identified as
a growth-fragmentation process whose Eve particle process is determined explicitly. 
Note that $\ell^x$ is a measurable function of $\mathscr{X}_x$: By \cite[Proposition 26]{LGR},
$\ell^x$ can be written, up to a multiplicative constant, as the limit of $\delta^{3/2}$ times the number
of components of $\mathscr{X}_x$ greater than $\delta$, when $\delta\to 0$. 
Similarly, Lemma \ref{ident-derivative} below shows that $\dot\ell^x$ is equal to twice 
the suitably renormalized sum of the components of $\mathscr{X}_x$
(some renormalization is needed because the sum is infinite). However, despite the fact that
$(\mathscr{X}_x)_{x\geq 0}$ is a Markov process with known distribution, it does not seem
easy to infer from this that the process $(\ell^x,\dot\ell^x)_{x\geq 0}$ is also Markov. 

The recent paper of Chapuy and Marckert \cite{CM} deals with the random measure ISE and addresses topics closely related
to the present work with very different (combinatorial) methods based on discrete approximations.
In particular, \cite{CM} proves that the density of ISE is continuously differentiable and discusses the regularity 
of the derivative. The study of discrete analogs also leads \cite{CM} to conjecture that
the derivative of the density satisfies a stochastic differential equation involving the density itself and
the distribution function of ISE (that is, the integral of the density over $(-\infty,t]$). One may
observe that conditioning the total volume of $\t$ to be equal to $1$ (as in the definition of ISE) makes it hopeless to
get a Markov property of the type of Theorem \ref{mainT}. 

The paper is organized as follows. Section 2 gives a number of preliminaries including a precise definition and
properties of the ``exit measure process'' $(\mathcal{X}^h_r)_{r\geq 0}$ and of the measures $\N^{*,z}_h$ that 
are used to describe the distribution of excursions
above or below the level $h$. In Section 3, we briefly recall the relations between super-Brownian motion and
our model of Brownian motion indexed by the Brownian tree, and we explain how Sugitani's results in \cite{Sug}
can be used to study the regularity of the process $(\ell^x)_{x\in\R}$ (more precise results about this
regularity are derived in \cite{CM}). Section 4 is devoted to technical estimates about the measures $\N^{*,z}_h$,
which play an important role in the subsequent proofs. The proof of Theorem \ref{mainT}
is given in Section 5. Section 6 then explains how the same method of proof can be
used to derive Theorem \ref{main-super}. Finally, Section 7 gives several open questions and complements.
In particular, we explain how Theorem \ref{mainT} provides information about the
model of random geometry known as the Brownian sphere, which has been studied
extensively in the recent years.

%

\section{Preliminaries}
\label{sec:preli}

\subsection{Snake trajectories}
\label{sna-tra}

We use the formalism of snake trajectories 
and we recall the main definitions that will be needed below.
We refer to \cite{ALG} for more information. 
A (one-dimensional) finite path $\w$ is a continuous mapping $\w:[0,\zeta]\la\R$, where the
number $\zeta=\zeta_{(\w)}$ is called the lifetime of $\w$. The space $\W$ of all finite paths is a Polish space when equipped with the
distance
$$d_\W(\w,\w')=|\zeta_{(\w)}-\zeta_{(\w')}|+\sup_{t\geq 0}|\w(t\wedge
\zeta_{(\w)})-\w'(t\wedge\zeta_{(\w')})|.$$
The endpoint or tip of the path $\w$ is denoted by $\wh \w=\w(\zeta_{(\w)})$.
For every $x\in\R$, we set $\W_x=\{\w\in\W:\w(0)=x\}$. The trivial element of $\W_x$ 
with zero lifetime is identified with the point $x$ of $\R$. 

\begin{definition}
\label{def:snakepaths}
Let $x\in \R$. A snake trajectory with initial point $x$ is a continuous mapping $s\mapsto \omega_s$
from $\R_+$ into $\W_x$ 
which satisfies the following two properties:
\begin{enumerate}
\item[\rm(i)] We have $\omega_0=x$ and the number $\sigma(\omega):=\sup\{s\geq 0: \omega_s\not =x\}$,
called the duration of the snake trajectory $\omega$,
is finite (by convention $\sigma(\omega)=0$ if $\omega_s=x$ for every $s\geq 0$). 
\item[\rm(ii)] For every $0\leq s\leq s'$, we have
$\omega_s(t)=\omega_{s'}(t)$ for every $t\in[0,\displaystyle{\min_{s\leq r\leq s'}} \zeta_{(\omega_r)}]$.
\end{enumerate} 
\end{definition}

Property (ii) is called the {\it snake property}.
We will write $\S_x$ for the set of all snake trajectories
with initial point $x$, and $\S$ for the union of the sets $\S_x$ for all $x\in\R$. If $\omega\in \S$, we often write $W_s(\omega)=\omega_s$ and $\zeta_s(\omega)=\zeta_{(\omega_s)}$
for every $s\geq 0$, and we omit $\omega$ in the notation. The sets $\S$ and $\S_x$ are equipped with the distance
$$d_{\S}(\omega,\omega')= |\sigma(\omega)-\sigma(\omega')|+ \sup_{s\geq 0} \,d_\W(W_s(\omega),W_{s}(\omega')).$$
A snake trajectory $\omega$ is completely determined 
by the knowledge of the lifetime function $s\mapsto \zeta_s(\omega)$ and of the tip function $s\mapsto \wh W_s(\omega)$: See \cite[Proposition 8]{ALG}.
For $\omega\in\S_x$ and $a\in\R$, we will use the obvious notation $\omega + a\in\S_{x+a}$.

Let $\omega\in \S$ be a snake trajectory and $\sigma=\sigma(\omega)$. The lifetime function $s\mapsto \zeta_s(\omega)$ codes a
compact $\R$-tree, which will be denoted 
by $\t=\t(\omega)$ and called the {\it genealogical tree} of the snake trajectory. This $\R$-tree is the quotient space $\t:=[0,\sigma]/\!\sim$ 
of the interval $[0,\sigma]$
for the equivalence relation
$$s\sim s'\ \hbox{if and only if }\ \zeta_s=\zeta_{s'}= \min_{s\wedge s'\leq r\leq s\vee s'} \zeta_r,$$
and $\t$ is equipped with the distance induced by
$$d_\zeta(s,s')= \zeta_s+\zeta_{s'}-2 \min_{s\wedge s'\leq r\leq s\vee s'} \zeta_r.$$
(notice that $d_\zeta(s,s')=0$ if and only if $s\sim s'$, and see e.g.~\cite{probasur} for more information about the
coding of $\R$-trees by continuous functions).  We let $p_{(\omega)}:[0,\sigma]\la \t$ stand
for the canonical projection. By convention, $\t$ is rooted at the point
$\rho:=p_{(\omega)}(0)=p_{(\omega)}(\sigma)$, and the volume measure $\mathrm{Vol}_{(\omega)}(\cdot)$ on $\t$ is defined as the pushforward of
Lebesgue measure on $[0,\sigma]$ under $p_{(\omega)}$ --- we often write $\mathrm{Vol}$ instead of $\mathrm{Vol}_{(\omega)}$.  As usual, for $a,b\in\t$, we say that $a$ is an ancestor
of $b$, or $b$ is a descendant of $a$, if $a$ belongs to the line segment from $\rho$ to $b$ in $\t$. 

By property (ii) in the definition of  a snake trajectory, the condition $p_{(\omega)}(s)=p_{(\omega)}(s')$ implies that 
$W_s(\omega)=W_{s'}(\omega)$. So the mapping $s\mapsto W_s(\omega)$ can be viewed as defined on the quotient space $\t$.
For $a\in \t$, we set $V_a(\omega):=\wh W_s(\omega)$ whenever  $s\in[0,\sigma]$ is such that $a=p_{(\omega)}(s)$ --- by the previous observation this does not
depend on the choice of $s$. We interpret $V_a$ as a ``label'' assigned to the ``vertex'' $a$ of $\t$. 
Notice that the mapping $a\mapsto V_a$ is continuous on $\t$. 
We will use the notation 
\begin{align*}
&W_*:=\min\{W_s(t): s\geq 0,t\in[0,\zeta_s]\}= \min\{V_a:a\in\t\},\\
&W^*:=\max\{W_s(t): s\geq 0,t\in[0,\zeta_s]\}= \max\{V_a:a\in\t\},
\end{align*}
and we also let $\mathcal{Y}(\omega)$ be the finite measure on $\R$ defined by setting
\begin{equation}
\label{occu-measure}
\langle\mathcal{Y},g\rangle = \int_0^{\sigma} g(\wh W_s)\,\dd s=\int_{\t} g(V_a)\,\mathrm{Vol}(\dd a),
\end{equation}
for any bounded continuous function $g:\R\la \R_+$. Trivially, $\mathcal{Y}$ is supported on $[W_*,W^*]$. 
%

\subsection{Re-rooting and truncation of snake trajectories}
\label{sec:reroot}

We now introduce two important operations on snake trajectories. The first one 
is the re-rooting operation (see \cite[Section 2.2]{ALG}). Let $\omega\in \S_x$ and
$r\in[0,\sigma(\omega)]$. Then $\omega^{[r]}$ is the snake trajectory in $\S_{\hat\omega_r}$ such that
$\sigma(\omega^{[r]})=\sigma(\omega)$ and for every $s\in [0,\sigma(\omega)]$,
\begin{align*}
\zeta_s(\omega^{[r]})&= d_\zeta(r,r\oplus s),\\
\wh W_s(\omega^{[r]})&= \wh W_{r\oplus s}(\omega),
\end{align*}
where we use the notation $r\oplus s=r+s$ if $r+s\leq \sigma$, and $r\oplus s=r+s-\sigma$ otherwise. 
By a remark following the definition of snake trajectories, these prescriptions completely determine $\omega^{[r]}$.

The genealogical tree $\t(\omega^{[r]})$ is then interpreted as the tree $\t(\omega)$ re-rooted at the vertex $p_{(\omega)}(r)$: More precisely,
the mapping $s\mapsto r\oplus s$ induces an isometry from $\t(\omega^{[r]})$
onto $\t(\omega)$, which maps the root of $\t(\omega^{[r]})$ to $p_{(\omega)}(r)$. Furthermore, the vertices
of $\t(\omega^{[r]})$ receive the same labels as in $\t(\omega)$.  

The second operation is the truncation of snake trajectories. For any $\w\in\W_x$ and $y\in\R$, we set
$$\tau_y(\w):=\inf\{t\in(0,\zeta_{(\w)}]: \w(t)=y\}\,,$$
with the usual convention $\inf\varnothing =\infty$.
Then if $\omega\in \S_x$ and $y\in \R$, we set, for every $s\geq 0$,
$$\nu_s(\omega):=\inf\Big\{t\geq 0:\int_0^t \mathrm{d}u\,\mathbf{1}_{\{\zeta_{(\omega_u)}\leq\tau_y(\omega_u)\}}>s\Big\}$$
(note that the condition $\zeta_{(\omega_u)}\leq\tau_y(\omega_u)$ holds if and only if $\tau_y(\omega_u)=\infty$ or $\tau_y(\omega_u)=\zeta_{(\omega_u)}$).
Then, setting $\omega'_s=\omega_{\nu_s(\omega)}$ for every $s\geq 0$ defines an element $\omega'$ of $\S_x$,
which will be denoted by  ${\rm tr}_y(\omega)$ and called the truncation of $\omega$ at $y$
(see \cite[Proposition 10]{ALG}). Informally, the truncated snake trajectory $\mathrm{tr}_y(\omega)$ behaves exactly like the
snake trajectory $\omega$ when the paths $\omega_s$ do not hit $y$, and the effect of the time 
change $\nu_s(\omega)$ is to ``eliminate'' those paths $\omega_s$ that hit $y$ and survive for a positive
amount of time.  The genealogical tree of 
${\rm tr}_y(\omega)$ is canonically and isometrically identified 
with the closed subset of $\t(\omega)$ consisting of all $a$ such that
$V_b(\omega)\not =y$ for every strict ancestor $b$ of $a$
(different from $\rho$ when $y=x$). 

Finally, for $\omega\in\S_x$ and $y\in\R\backslash\{x\}$, we  define the excursions of $\omega$ away from $y$. We let $(\alpha_j,\beta_j)$, $j\in J$, be the connected components of the open set
$$\{s\in[0,\sigma]:\tau_y(\omega_s)<\zeta_{(\omega_s)}\}$$
(note that the indexing set $J$ may be empty).
We notice that $\omega_{\alpha_j}=\omega_{\beta_j}$ for every $j\in J$, by the snake property, and $\wh\omega_{\alpha_j}=y$.
For every $j\in J$, we define a snake trajectory $\omega^j\in\S_y$ by setting
$$\omega^j_{s}(t):=\omega_{(\alpha_j+s)\wedge\beta_j}(\zeta_{(\omega_{\alpha_j})}+t)\;,\hbox{ for }0\leq t\leq \zeta_{(\omega^j_s)}
:=\zeta_{(\omega_{(\alpha_j+s)\wedge\beta_j})}-\zeta_{(\omega_{\alpha_j})}\hbox{ and } s\geq 0.$$
We say that $\omega^j$, $j\in J$, are the excursions of $\omega$ away from $y$. 

\subsection{The Brownian snake excursion measure}
\label{sna-mea}

Let $x\in\R$. The Brownian snake excursion 
measure $\N_x$ is the $\sigma$-finite measure on $\S_x$ that satisfies the following two properties: Under $\N_x$,
\begin{enumerate}
\item[(i)] the distribution of the lifetime function $(\zeta_s)_{s\geq 0}$ is the It\^o 
measure of positive excursions of linear Brownian motion, normalized so that, for every $\ve>0$,
$$\N_x\Big(\sup_{s\geq 0} \zeta_s >\ve\Big)=\frac{1}{2\ve};$$
\item[(ii)] conditionally on $(\zeta_s)_{s\geq 0}$, the tip function $(\wh W_s)_{s\geq 0}$ is
a Gaussian process with mean $x$ and covariance function 
$$K(s,s')= \min_{s\wedge s'\leq r\leq s\vee s'} \zeta_r.$$
\end{enumerate}
Informally, the lifetime process $(\zeta_s)_{s\geq 0}$ evolves under $\N_x$ like a Brownian excursion,
and conditionally on $(\zeta_s)_{s\geq 0}$, each path $W_r$ is a linear Brownian path started from $x$ with
lifetime $\zeta_r$, which
is ``erased'' from its tip when $\zeta_r$ decreases, and is ``extended'' when $\zeta_r$ increases.
The measure $\N_x$ can be interpreted as the excursion measure away from $x$ for the 
Markov process in $\W_x$ called the (one-dimensional) Brownian snake. Note that the preceding informal description 
applies as well to the Brownian snake, except that, in that case, the lifetime process evolves like
a reflecting Brownian motion in $[0,\infty)$.
We refer to 
\cite{Zurich} for a detailed study of the Brownian snake with a more general
underlying spatial motion. 

As usual for excursion measures,
we can state a Markov property under $\N_x$. Let $u>0$ and let $F$ and $H$ be two nonnegative
measurable functions defined respectively on the space of all continuous functions from $[0,u]$
into $\W_x$ and on the space of all continuous functions from $[0,\infty)$
into $\W_x$. Then,
\begin{align}
\label{Markov-excu}
&\N_x\Big(\mathbf{1}_{\{u<\sigma\}}\,F\big((W_r)_{0\leq r\leq u}\big)\,H\big((W_{u+s})_{s\geq 0}\big)\Big)\\
&\qquad= \N_x\Big(\mathbf{1}_{\{u<\sigma\}}\,F\big((W_r)_{0\leq r\leq u}\big)\,\E^*_{W_u}\Big[H\big((W_{s})_{s\geq 0}\big)\Big]\Big),\nonumber
\end{align}
where, for every $\w\in\W_x$, $\P^*_\w$ denotes the law of the Brownian snake started from $\w$
and stopped when the lifetime process hits $0$ (see \cite[Section IV.4]{Zurich}).

For every $r>0$, we have
$$\N_x(W^*>x+r)=\N_x(W_*<x-r)={\displaystyle \frac{3}{2r^2}}$$ (see e.g. \cite[Section VI.1]{Zurich}). 
In particular, $\N_x(y\in[W_*,W^*])<\infty$ if $y\not =x$. 
We will use the first-moment formula under $\N_x$, which states that, for any
nonnegative measurable function $F$ on $\mathcal{W}_x$, 
\begin{equation}
\label{first-mom}
\N_x\Big(\int_0^\sigma F(W_s)\,\dd s\Big) = \int_0^\infty \dd t\,\mathbf{E}_x\big[ F\big((B_r)_{0\leq r\leq t}\big)\big],
\end{equation}
where $B$ denotes a linear Brownian motion that starts from $x$ under the probability measure $\mathbf{P}_x$
(see \cite[Chapter 4]{Zurich}). We also recall the re-rooting invariance property of $\N_0$ \cite[Theorem 2.3]{LGW}. 
To state this property, it is convenient to modify a little the definition of a re-rooted snake trajectory
in the preceding section: if $\omega\in\S_0$ and $r\in[0,\sigma(\omega)]$, we set
$\wt\omega^{[r]}=\omega^{[r]}-\wh \omega_r$ (we just shift the snake trajectory $\omega^{[r]}$
so that it belongs to $\S_0$ instead of $\S_{\hat\omega_r}$). Then, for any
nonnegative measurable function $F$ on $[0,\infty)\times \S_0$, we have
\begin{equation}
\label{re-root-for}
\N_0\Bigg(\int_0^\sigma \dd r\,F(r,\wt\omega^{[r]})\Bigg)= \N_0\Bigg(\int_0^\sigma \dd r\,F(r,\omega)\Bigg).
\end{equation}

The following scaling property is often useful. For $\lambda>0$, for every 
$\omega\in \S_x$, we define $\theta_\lambda(\omega)\in \S_{x\sqrt{\lambda}}$
by $\theta_\lambda(\omega)=\omega'$, with
$$\omega'_s(t):= \sqrt{\lambda}\,\,\omega_{s/\lambda^2}(t/\lambda)\,,\ 
\hbox{for } s\geq 0\hbox{ and }0\leq t\leq \zeta'_s:=\lambda\zeta_{s/\lambda^2}.$$
Then $\theta_\lambda(\N_x)= \lambda\, \N_{x\sqrt{\lambda}}$. 

Let us now define exit measures. We argue under $\N_x$, and fix $y\in\R\backslash\{x\}$. Then, the idea is to make sense of a quantity that ``measures'' the number of paths 
$W_s$ that hit level $y$ and are stopped at that hitting time. Precisely, one shows \cite[Proposition 34]{Disks} that the limit
\begin{equation}
\label{formu-exit}
\mathcal{L}^y_t:=\lim_{\ve \to 0} \frac{1}{\ve^2} \int_0^t \dd s\,\mathbf{1}_{\{\zeta_s\leq\tau_y(W_s),\,|\wh W_s-y|<\ve\}}
\end{equation}
exists for every $t\in[0,\sigma]$, $\N_x$ a.e., and defines a continuous increasing function
called the exit local time from $(y,\infty)$ (if $x>y$) or from $(-\infty,y)$ (if $y>x$). 
The exit measure is then defined by $\z_y:=\mathcal{L}^y_\sigma$, and we have  $\z_y>0$ if and only if $y\in[W_*,W^*]$, $\N_x$ a.e. 
This definition of the exit local time and of $\z_y$
is a particular case of the theory of exit measures, see \cite[Chapter V]{Zurich} where a different but equivalent approximation of
$\mathcal{L}^y_t$ is used. It follows from the approximation \eqref{formu-exit} that $\z_y$ is $\N_x$ a.e.~equal 
to a measurable function of the truncated snake $\mathrm{tr}_y(\omega)$. We will use the following formula, for every $\lambda> 0$,
\begin{equation}
\label{Laplace-exit}
\N_x\Big(1-\exp (-\lambda \z_y)\Big)= \Big(|x-y|\sqrt{2/3} + \lambda^{-1/2}\Big)^{-2}.
\end{equation}
See formula (6) in \cite{CLG} for a brief justification. In particular, we have $\N_x(\z_y)=1$. 

%
We now recall the special Markov property of the Brownian snake under $\N_0$ (see in particular the appendix
of \cite{subor}).  

\begin{proposition}[Special Markov property]
\label{SMP}
Let $x\in \R$ and $y\in\R\backslash\{x\}$. Under the measure $\N_x(\dd \omega)$, 
let $\omega^j$, $j\in J$, be the excursions of $\omega$ away from $y$ and consider the point measure
$$\n_y=\sum_{j\in J} \delta_{\omega^j}.$$
Then, under the probability measure $\N_x(\dd\omega\,|\, y\in[W_*,W^*])$ and conditionally on $\z_y$, the point measure
$\n_y$ is Poisson with intensity $\z_y\,\N_y(\cdot)$ and is independent of $\mathrm{tr}_y(\omega)$.
\end{proposition}

\subsection{The exit measure process at a point}
\label{sec:exitlocal}

Let us fix a point $x\in\R$ in this section and consider a snake trajectory $\omega$ distributed according to $\N_x$. 
An important role in this work will be played by a process $(\mathcal{X}^x_r)_{r>0}$, such that
for every $r>0$, $\mathcal{X}^x_r$ measures the
``quantity'' of paths $W_s(\omega)$ that have accumulated a
local time at $x$ exactly equal to $r$. 
The precise 
definition of $\mathcal{X}^x_r$ belongs to
the general theory of exit measures
and we refer to the introduction of \cite{ALG}
for more details (roughly speaking, one needs to consider the Brownian snake 
whose spatial motion is the pair consisting of a
linear Brownian motion and its local time at $x$,
and then the exit measure from the set $\R\times [0,r)$). 
One proves that the process $(\mathcal{X}^x_r)_{r> 0}$
is distributed under $\N_x$ according to the excursion measure 
of the continuous-state branching process with branching
mechanism $\psi(u)=\sqrt{8/3}\,u^{3/2}$ (in short, the $\psi$-CSBP, we refer
to \cite[Chapter II]{Zurich}
for basic facts about CSBPs).
This means that, if $\n=\sum_{k\in K}\delta_{\omega_k}$
is a Poisson point measure with intensity $\alpha\,\N_x$, the
process $X$ defined by $X_0=\alpha$ and, for every $r>0$,
$$X_r:=\sum_{k\in K} \mathcal{X}^x_r(\omega_k),$$
is a $\psi$-CSBP started at $\alpha$ (see \cite[Section 2.4]{LGR}). In particular,
$(\mathcal{X}^x_r)_{r> 0}$
has a c\`adl\`ag modification under $\N_x$, which
we consider from now on. We take $\mathcal{X}^x_0=0$ 
by convention and call $(\mathcal{X}^x_r)_{r\geq 0}$ the 
exit measure process at $x$.

Still under $\N_x$, we can also define the exit measure process at $y$
for any $y\not=x$. We can either rely on the general theory
of exit measures, or use the point process $\sum_{j\in J} \delta_{\omega^j}$
of excursions away from $y$ (as in Proposition \ref{SMP}) to define
for every $r>0$,
$$\mathcal{X}^y_r:=\sum_{j\in J} \mathcal{X}^y_r(\omega^j)$$
(note that the quantities $\mathcal{X}^y_r(\omega^j)$ make sense 
by the special case $y=x$ treated before). 
We also set $\mathcal{X}^y_0=\z_y$. It follows from 
Proposition \ref{SMP} 
and the preceding considerations that, 
under the probability measure $\N_x(\cdot\mid y\in[W_*,W^*])=\N_x(\cdot\mid \z_y>0)$, conditionally on $\z_y$,
the process $(\mathcal{X}^y_r)_{r\geq 0}$ is a 
$\psi$-CSBP started at $\z_y$ and is independent of $\tr_y(\omega)$. 
Again we call $(\mathcal{X}^y_r)_{r\geq 0}$ the 
exit measure process at $y$.

\subsection{The positive excursion measure}

Under $\N_0$, the paths $\omega_s$, $0<s<\sigma$, take both positive and negative values, simply because they
behave like one-dimensional Brownian paths started from $0$. We will now introduce another important measure
on $\S_0$, which is supported on snake trajectories taking only nonnegative values. For $\delta\geq 0$, let $\S_0^{(\delta)}$
be the set of all $\omega\in \S_0$ such that 
$\sup_{s\geq 0}(\sup_{t\in[0,\zeta_s(\omega)]}|\omega_s(t)|)> \delta$. 
Also set
$$\S_0^+=\{\omega\in\S_0: \omega_s(t)\geq 0 \hbox{ for every }s\geq 0\hbox{ and }t\in[0,\zeta_s(\omega)]\}\cap \S^{(0)}_0.$$
There exists a $\sigma$-finite measure $\N^*_0$ on
$\S_0$, which is supported on  $\S_0^+$, 
and gives finite mass to the sets $\S_0^{(\delta)}$ for every $\delta>0$, such that
$$\N^*_0(G)=\lim_{\ve\to 0}\frac{1}{\ve}\,\N_\ve(G(\mathrm{tr}_{0}(\omega))),$$
for every bounded continuous function $G$ on $\S_0$ that vanishes 
on $\S_0\backslash\S_0^{(\delta)}$ for some $\delta>0$ (see \cite[Theorem 23]{ALG}). 
Under $\N^*_0$, each of the paths $\omega_s$, $0<s<\sigma$, starts from $0$, then
stays positive during some time interval $(0,\alpha)$, and is stopped immediately
when it returns to $0$, if it does return to $0$. 

In a way analogous to the definition of exit measures, one can make sense of the ``quantity'' of paths $\omega_s$ that return to $0$ under $\N^*_0$.
To this end, one proves that the limit
\begin{equation}
\label{approxz*0}
\z^*_0:=\lim_{\ve\to 0} \frac{1}{\ve^2} \int_0^\sigma \dd s\,\mathbf{1}_{\{0<\wh W_s<\ve\}}
\end{equation}
exists $\N^*_0$ a.e. See \cite[Proposition 30]{ALG} for a slightly weaker result --- the stronger form
stated above follows from the results of \cite[Section 10]{Disks}. 
Notice that replacing the 
limit by a liminf in \eqref{approxz*0} allows us to make sense of $\z^*_0(\omega)$ for every $\omega\in\S_0^+$.

The following conditional versions of the
measure $\N^*_0$ play a fundamental role in the present work. According to
\cite[Proposition 33]{ALG}, there exists a unique collection $(\N^{*,z}_0)_{z>0}$ of probability measures on $\S_0^+$
such that:
\begin{enumerate}
\item[\rm(i)] We have
$\;\N^*_0={\displaystyle
 \sqrt{\frac{3}{2 \pi}} 
\int_0^\infty \mathrm{d}z\,z^{-5/2}\, \N^{*,z}_0}$.
\item[\rm(ii)] For every $z>0$, $\N^{*,z}_0$ is supported on $\{\z^*_0=z\}$.
\item[\rm(iii)] For every $z,z'>0$, $\N^{*,z'}_0=\theta_{z'/z}(\N^{*,z}_0)$.
\end{enumerate}
Informally, $\N^{*,z}_0=\N^*_0(\cdot\mid \z^*_0=z)$. It will be convenient to
write
\begin{equation}
\label{Levy-mea}
\bn(\dd z)=\sqrt{\frac{3}{2 \pi}} 
\,z^{-5/2}\,\mathrm{d}z,
\end{equation}
so that $\bn(\dd z)$ is the ``law'' of $\z^*_0$ under $\N^*_0$.
We note that the convergence \eqref{approxz*0}
also holds $\N^{*,z}_0$ a.s., with $\z^*_0$ replaced by $z$ (\cite[Corollary 37]{Disks}), and we record the formula
\begin{equation}
\label{N*sigma}
\N^{*,z}_0(\sigma)=z^2,
\end{equation}
for every 
$z>0$ (see e.g.~\cite[Proposition 10]{LGR}). 

It will be convenient to write $\check\N^{*,z}_0$
for the pushforward of $\N^{*,z}_0$ under the mapping $\omega\to -\omega$. Furthermore, for
every $h\in\R$, we write $\N^{*,z}_h$, resp. $\check\N^{*,z}_h$ for the pushforward of $\N^{*,z}_0$,
resp. of $\check\N^{*,z}_0$, under the shift $\omega\mapsto \omega + h$. 

The next theorem relates the measures $\N_x$ and $\N^*_0$ via a re-rooting transformation.
Recall that, for every $\omega\in \S$ and every $s\in[0,\sigma(\omega)]$, 
$\omega^{[s]}$ denotes the snake trajectory $\omega$ re-rooted at $s$ (Section \ref{sec:reroot}).

\begin{theorem}{\rm\cite[Theorem 28]{ALG}}
\label{re-root-theo}
Let $G$ be a nonnegative measurable function on $\S$. Then,
$$\N^*_0\Bigg(\int_0^\sigma \dd r\, G(\omega^{[r]})\Bigg)
= 2\int_0^\infty \dd b\,\N_b\Big( \z_0\,G(\tr_0(\omega))\Big).
$$
\end{theorem}

As a first application, we can take $G(\omega)=g(\omega(0))$ where 
$g:\R\la\R_+$ is measurable. Since $\N_b(\z_0)=1$ for every $b>0$, it follows that
\begin{equation}
\label{first-mo}
\N^*_0\Bigg(\int_0^\sigma \dd r\, g(\wh W_r)\Bigg)= 2\int_0^\infty \dd b\,g(b).
\end{equation}

\subsection{Excursion theory}
\label{sec:excu}

We now recall the main theorem of the excursion theory developed 
in \cite{ALG}. We fix $x\in\R$ and $y\in\R$. 
We consider a random snake trajectory $\omega$ distributed according to $\N_x$.
The goal of this excursion theory is to describe 
the connected components of $\{v\in\t(\omega):V_v\not =y\}$, and the 
evolution of labels on these connected components (there is an obvious 
analogy with classical excursion theory for linear Brownian motion). 
Let $\mathcal{C}$ be a connected component of $\{v\in\t(\omega):V_v\not =y\}$,
and exclude the case where $\mathcal{C}$ contains the root $\rho$ of $\t(\omega)$
(this case occurs when $y\not=x$).
If $\ov{\mathcal{C}}$ denotes the closure of $\mathcal{C}$, there
is a unique point $u$ of $\ov{\mathcal{C}}$ at minimal distance from the root
(in such a way that all points of $\ov{\mathcal{C}}$ are descendants of $u$)
and we have $V_u=y$.
Following \cite{ALG}, we say that $u$ is an excursion debut from $y$. 
We observe that a branching point of $\t(\omega)$ cannot be an excursion debut from $y$: indeed, branching points of $\t(\omega)$ correspond under 
$p_{(\omega)}$ to times of local minimum of $(0,\sigma)\ni s\mapsto \zeta_s$, and, if $r$ is such a time of local minimum, the 
conditional distribution of $(\wh W_s)_{s\geq 0}$ given $(\zeta_s)_{s\geq 0}$ shows that $\wh W_r\not = y$ a.s., which gives
the desired result since there are only countably many times of local minimum.

Let $\mathcal{C}$ be a connected component of $\{v\in\t(\omega):V_v\not =y\}$ and let $u$ be
the associated excursion debut. We can 
code $\mathcal{C}$ and the labels on $\mathcal{C}$
via a snake trajectory which is defined as follows. First we observe that there are exactly
 two times $s_0<s'_0$
such that $p_{(\omega)}(s_0)=p_{(\omega)}(s'_0)=u$,
and the set 
$p_{(\omega)}([s_0,s'_0])$ is the subtree of all descendants 
of $u$ (here we use the
fact that $u$ is not a branching point).  We then define a snake trajectory $\tilde\omega^{(u)}\in \S_0$
coding the subtree $p_{(\omega)}([s_0,s'_0])$ (and its labels) by setting
$$\tilde\omega^{(u)}_s(t):=\omega_{(s_0+s)\wedge s'_0}(\zeta_{s_0}+t)\,\hbox{ for }0\leq t\leq \zeta_{(s_0+s)\wedge s'_0}.$$
We finally set $\omega^{(u)}:=\tr_y(\tilde\omega^{(u)})$
and we observe that the compact $\R$-tree $\ov{\mathcal{C}}$ is identified 
isometrically to the tree $\t(\omega^{(u)})$, and moreover this identification
preserves the labels. Also, the restriction of the volume measure $\mathrm{Vol}_{(\omega)}$ to $\mathcal{C}$ (or to $\ov{\mathcal{C}}$)
corresponds via the latter identification to the volume measure $\mathrm{Vol}_{(\omega^{(u)})}$ of $\t(\omega^{(u)})=\ov{\mathcal{C}}$.
Note that $\mathrm{Vol}_{(\omega)}$ puts no mass on $\ov{\mathcal{C}}\backslash\mathcal{C}$ since 
$\ov{\mathcal{C}}\backslash\mathcal{C}\subset\{v\in\t(\omega): V_v=y\}$.

We say that
$\omega^{(u)}$ is an excursion above $y$ if the values of $V_v$ for $v\in\mathcal{C}$ are greater than $y$
(equivalently the paths $\omega^{(u)}_s$ take values in $[y,\infty)$), and  that $\omega^{(u)}$ is an excursion
below  $y$  if the values of $V_v$ for $v\in\mathcal{C}$ are smaller than $y$. We note that the terminology is a bit misleading, since an excursion away from $y$,
as considered in Proposition \ref{SMP}, will contain infinitely many
excursions above or below $y$.

Recall from Section \ref{sec:exitlocal} 
the definition of the exit measure process at $y$, which is denoted by $(\mathcal{X}^y_r)_{r\geq 0}$.
If $y\notin[W_*,W^*]$ (which does not occur when $y=x$, and is equivalent to $\z_y=0$ when $y\not=x$),  
there are no excursion debuts from $y$. For this reason, we suppose that $\z_y>0$ when $y\not=x$  in the following
lines. By Proposition 3 of \cite{ALG} (and an application of the special Markov property when $y\not=x$), excursion debuts 
from $y$ are in one-to-one 
correspondence with the jump times of the process $(\mathcal{X}^y_r)_{r\geq 0}$,
in such a way that, if $u$ is an excursion debut and $s\in[0,\sigma]$ is such that
$p_{(\omega)}(s)=u$, the associated jump time of the exit measure process at $y$
is the total local time at $y$ accumulated by the path $\omega_s$. We can 
list the jump times of $(\mathcal{X}^y_r)_{r\geq 0}$ in a sequence $(r_i)_{i\in\N}$
in decreasing order of the jumps $\Delta\mathcal{X}^y_{r_i}=\mathcal{X}^y_{r_i}-\mathcal{X}^y_{r_i-}$.
For every $i\in\N$, we write $u_i$ for the excursion debut associated with $r_i$.

The following theorem is essentially Theorem 4 of \cite{ALG}. We write $\N_x^{(y)}=\N_x(\cdot \mid \z_y>0)$ when $y\not=x$,
and $\N_x^{(x)}=\N_x$.

\begin{theorem}
\label{theo-excursion}
Under $\N_x^{(y)}$, conditionally on $(\mathcal{X}^y_r)_{r\geq 0}$,
the excursions $ \omega^{(u_i)}$, $i\in \N$, are independent, and independent
of $\tr_y(\omega)$, and the 
conditional distribution of $\omega^{(u_i)}$ is 
$$\frac{1}{2}\Big(\N^{*,\Delta\mathcal{X}^y_{r_i}}_y+\check\N^{*,\Delta\mathcal{X}^y_{r_i}}_y\Big),$$
where $\Delta\mathcal{X}^y_{r_i}=\mathcal{X}^y_{r_i}-\mathcal{X}^y_{r_i-}$ is the jump of $\mathcal{X}^y$ at time $r_i$.
\end{theorem}

To be specific, Theorem 4 of \cite{ALG} deals with the case $y=x$ (in that case, $\tr_y(\omega)$
is trivial), but
then an application of the special Markov property (Proposition \ref{SMP}) yields the case $y\not=x$. 

\subsection{A path transformation of L\'evy processes}
\label{path-transfo}

The classical Lamperti transformation \cite{Lam} shows that the continuous-state branching process
$(\XX^y_r)_{r\geq 0}$ of the preceding section can be obtained as a time change
of a stable L\'evy process with no negative jumps. In this section, we state a path transformation 
of L\'evy processes that will be relevant in forthcoming proofs. Let $\beta\in(1,2)$,
and let $(U_s)_{s\geq 0}$ be a (centered) stable L\'evy process with index $\beta$ and with
no negative jumps, such that $U_0=a>0$. Then the Laplace transform of $U_s-a$ is well defined and
given by $\E[\exp(-\lambda (U_s-a))]=\exp(c\,s\,\lambda^\beta)$ for every $\lambda>0$, where $c>0$
is a constant. We say that the Laplace exponent of $U$ is $c\,\lambda^\beta$. For every $t>0$, we write $(U^{\mathrm{br},a,t}_s)_{0\leq s\leq t}$ for the 
associated bridge of duration $t$ from $a$ to $0$, that is, for the process $(U_s)_{0\leq s\leq t}$ 
conditioned on $U_t=0$. We refer to \cite{FPY} for a precise definition and construction
of this bridge. 

We then set
$$T_0:=\inf\{s\geq 0: U_s=0\},$$
and we consider the following transformation of the path of $U$
over the time interval $[0,T_0]$. Let $R$ be a nonnegative random variable
which is uniformly distributed over $[0,T_0]$ conditionally given the process $U$. For every $s\in[0,T_0]$, we set
$$\wt U_s=\left\{\begin{array}{ll} U_{R+s}-U_R+a\quad&\hbox{if }0\leq s\leq T_0-R;\\
U_{R+s-T_0}-U_R&\hbox{if }T_0-R\leq s\leq T_0.
\end{array}
\right.
$$

\begin{lemma}
\label{Kemp}
The conditional distribution of $(\wt U_s)_{0\leq s\leq T_0}$ knowing that $T_0=t$
is the law of $(U^{\mathrm{br},a,t}_s)_{0\leq s\leq t}$.
\end{lemma}

A  discrete version of the previous statement, for centered random walks 
with negative jumps of size $-1$ only, is easy to prove from the arguments based on the 
cyclic lemma that lead to the classical Kemperman lemma --- see e.g. Section 6.1 of \cite{Pit}.
Then Lemma \ref{Kemp} follows by applying a suitable invariance principle. Alternatively, Corollary 8 of \cite{BCP}
 gives the analog of Lemma \ref{Kemp} when $U$ is replaced by
a linear Brownian motion, and Section 4 of the same paper explains how this can 
be extended to the setting of processes with cyclically exchangeable increments (which 
is more than we need here). 

In order to apply Lemma \ref{Kemp}, we note that the collection of jumps of the 
process $\wt U$ over the time interval $[0,T_0]$ is the same as the collection 
of jumps of $U$ over the same interval. Write $\D(\R_+,\R)$ for the Skorokhod 
space of c\`adl\`ag real functions on $\R_+$, and $\D_0(\R_+,\R)$ for the subset of $\D(\R_+,\R)$
consisting of functions with compact support.
Then, if $g$ is a nonnegative measurable
function on $\R_+$, and $F$ is a nonnegative measurable function defined 
on $\D_0(\R_+,\R)$ such that $F(\w)$ only depends 
on the sequence of jumps of $\w$ ordered in nonincreasing size, we have
\begin{equation}
\label{Kemp-appli}
\E\Big[g(T_0)\,F( (U_{s\wedge T_0})_{s\geq 0})\Big]=\E\Big[g(T_0)\,F( (\wt U_{s\wedge T_0})_{s\geq 0})\Big]
=\int \pi_a(\dd t)\,g(t)\,\E\Big[F((U^{\mathrm{br},a,t}_{s\wedge t})_{s\geq 0})\Big],
\end{equation}
where $\pi_a$ stands for the law of $T_0$. In other words, the conditional distribution of the
sequence of jumps of $(U_s)_{0\leq s\leq T_0}$ (ordered in nonincreasing size) knowing that $T_0=t$
is the distribution of the sequence of jumps of $(U^{\mathrm{br},a,t}_s)_{0\leq s\leq t}$. 

\section{The connection with super-Brownian motion}
\label{sec:super}

In this section, we briefly recall the connection between the Brownian snake excursion measures $\N_x$
and super-Brownian motion, referring to \cite{Zurich} for more details. We fix $\alpha>0$, and consider a Poisson point measure on $\S$,
$$\n=\sum_{k\in K} \delta_{\omega_k}$$
with intensity $\alpha\,\N_0$. Then one can construct a 
one-dimensional super-Brownian motion $(\mathbf{X}_t)_{t\geq 0}$ 
with branching mechanism $\phi(u)=2u^2$ and initial value $\mathbf{X}_0=\alpha\delta_0$, such that,
for any nonnegative measurable function $g$ on $\R$,
\begin{equation}
\label{occu-super}
\int_0^\infty \langle \mathbf{X}_t,g\rangle \dd t= \sum_{k\in K} \langle \mathcal{Y}(\omega_k),g\rangle,
\end{equation}
where $\mathcal{Y}(\omega_k)$ is defined in formula \eqref{occu-measure}. 
In a more precise way, the process $(\mathbf{X}_t)_{t\geq 0}$ is defined by setting, for
every $t> 0$ and every nonnegative measurable function $g$ on $\R$,
$$\langle \mathbf{X}_t,g\rangle := \sum_{k\in K}  \int_0^{\sigma(\omega_k)} \,\dd_rl^t_r(\omega_k)\,g(\wh W_r(\omega_k)),$$
where $l^t_r(\omega_k)$ denotes the local time of the process $s\mapsto \zeta_s(\omega_k)$ at level $t$ and at time $r$, 
and the notation $\dd_rl^t_r(\omega_k)$ refers to integration with respect to the nondecreasing function $r\mapsto l^t_r(\omega_k)$
 (see Chapter 4 of \cite{Zurich}).
We are primarily interested in the total occupation measure
$$\mathbf{Y}:=\int_0^\infty \mathbf{X}_t \,\dd t.$$
It follows from the results of Sugitani \cite[Theorem 4]{Sug} that $\mathbf{Y}$ has (a.s.) a 
continuous density $(L^x)_{x\in\R}$ with respect to Lebesgue measure, and 
this density is continuously differentiable on $(0,\infty)$
and on $(-\infty,0)$. On the other hand, for every $\ve>0$, the event $A$ where the point
measure $\n$ has exactly one atom $\omega_{*}$ such that $W^*(\omega_*)\geq \ve$
has positive probability, and, conditionally on this event, $\omega_*$ is distributed
according to $\N_0(\cdot\mid W^*\geq\ve)$. Furthermore, on the event $A$, formula 
\eqref{occu-super} entails that the restriction of $\mathbf{Y}$ to $(\ve,\infty)$ coincides
with the restriction of $\mathcal{Y}(\omega_*)$ to the same set. It follows 
that, a.s. under the probability measure $\N_0(\cdot\mid W^*\geq\ve)$, $\mathcal{Y}$
has a continuously differentiable density on $(\ve,\infty)$. Since $\ve$ was arbitrary, 
and using a symmetry argument,
we easily conclude that $\mathcal{Y}$
has a continuously differentiable density on $(-\infty,0)\cup(0,\infty)$, $\N_0$ a.e.

In fact, we can remove the ``singularity'' at $0$.
Indeed, we may use the  re-rooting invariance property of $\N_0$ (formula \eqref{re-root-for}) to
obtain that $\mathcal{Y}$ has a continuously differentiable density on $\R\backslash\{x\}$, $\mathcal{Y}(\dd x)$ a.e., $\N_0$ a.e.
It follows that $\mathcal{Y}$
has a continuously differentiable density on $\R$, $\N_0$ a.e. --- as already mentioned, this fact also follows from the results of \cite{CM},
which are proved via a completely different method.
As in the introduction above, we write $(\ell^x,x\in\R)$ for the density of $\y$ (under $\N_0$)
and call $\ell^x$ the local time at level $x$. The derivative of $\ell^x$ is denoted by $\dot\ell^x$.

\section{Technical estimates}

The following lemma is a key ingredient of the proof of our main result.

\begin{lemma}
\label{tech-lemma}
{\rm(i)} For every $z>0$ and $\ve>0$,
$$\N^{*,z}_0\Bigg(\int_0^\sigma \dd s\,\mathbf{1}_{\{\wh W_s<\ve\}}\Bigg)= \ve^4\,f(\frac{z}{\ve^2}),$$
where the function $f:(0,\infty)\la(0,\infty)$ is continuous, and  $u^{-1}f(u)\la 1$ as $u\to\infty$.

\smallskip
\noindent {\rm(ii)} There exists a constant $C$ such that, for every $\alpha\in(0,1]$ and $\ve\in(0,\sqrt{\alpha}\,]$, we have
$$\ve^{-4}  \N^{*}_0\Bigg(\Big(\int_0^\sigma \dd s\,\mathbf{1}_{\{\wh W_s<\ve\}} \Big)^2
\,\mathbf{1}_{\{\z^*_0\leq \alpha\}}
\Bigg) \leq C\sqrt{\alpha}.$$
\end{lemma}

\proof (i) For every $z>0$ and $x>0$, set
$$p_z(x)= 2\,\Big(\frac{3}{2}\Big)^{3/2}\,\sqrt{\pi}\,z^{3/2}x^{-4}\,\Upsilon(\frac{3z}{2x^2}\Big),$$
where the function $\Upsilon$ is defined on $(0,\infty)$ by
$$\Upsilon(x)= \frac{2}{\sqrt{\pi}}(x^{1/2}+x^{-1/2}) -2(x+\frac{3}{2})e^x\,\mathrm{erfc}(\sqrt{x}),$$
with the usual notation $\mathrm{erfc}(y)=\frac{2}{\sqrt{\pi}}\int_y^\infty e^{-x^2}\dd x$. 
Note that $\Upsilon(x)=\frac{2}{\sqrt{\pi}}x^{-1/2} + O(1)$ as $x\to 0$, and $\Upsilon(x)=\frac{3}{2\sqrt{\pi}}x^{-5/2} + O(x^{-7/2})$ as $x\to\infty$.

By \cite[Proposition 3]{Spine}, we have for every $\ve>0$,
$$\N^{*,z}_0\Bigg(\int_0^\sigma \dd s\,\mathbf{1}_{\{\wh W_s<\ve\}}\Bigg)= 
z^2 \int_0^\ve \dd x\,p_z(x).$$
Using the preceding formula for $p_z(x)$, and the change of variable $y=\frac{3z}{2 x^2}$, we obtain
$$\N^{*,z}_0\Bigg(\int_0^\sigma \dd s\,\mathbf{1}_{\{\wh W_s<\ve\}}\Bigg)
= \sqrt{\pi}\,z^2\int_{(3z)/(2\ve^2)}^\infty \dd y\,\sqrt{y}\,\Upsilon(y).$$
We thus get the formula of part (i) of the lemma with the function $f$ defined by
$$f(u)=\sqrt{\pi} \,u^2\int_{3u/2}^\infty \dd y\,\sqrt{y}\,\Upsilon(y).$$
The properties of $f$ stated in the lemma follow from this explicit expression and the asymptotics of $\Upsilon(x)$ as $x\to\infty$.

\medskip
\noindent (ii) We use a scaling argument to write
\begin{align*}
&\ve^{-4}  \N^{*}_0\Bigg(\Big(\int_0^\sigma \dd s\,\mathbf{1}_{\{\wh W_s<\ve\}} \Big)^2
\,\mathbf{1}_{\{\z^*_0\leq \alpha\}}
\Bigg)\\
&\quad=\sqrt{\frac{3}{2\pi}}\, \ve^{-4}\int_0^\alpha \dd z\,z^{-5/2}\,\N^{*,z}_0\Bigg(\Big(\int_0^\sigma \dd s\,\mathbf{1}_{\{\wh W_s<\ve\}} \Big)^2\Bigg)\\
&\quad=\sqrt{\frac{3}{2\pi}}\, \ve^{4}\int_0^\alpha \dd z\,z^{-5/2}\,\N^{*,\ve^{-2}z}_0\Bigg(\Big(\int_0^\sigma \dd s\,\mathbf{1}_{\{\wh W_s<1\}} \Big)^2\Bigg)\\
&\quad=\sqrt{\frac{3}{2\pi}}\, \ve\int_0^{\ve^{-2}\alpha} \dd z\,z^{-5/2}\,\N^{*,z}_0\Bigg(\Big(\int_0^\sigma \dd s\,\mathbf{1}_{\{\wh W_s<1\}} \Big)^2\Bigg)\\
&\quad=\ve J(\frac{\alpha}{\ve^2}),
\end{align*}
where we have set, for every $a>0$, 
$$J(a):=\N^{*}_0\Bigg(\Big(\int_0^\sigma \dd s\,\mathbf{1}_{\{\wh W_s<1\}} \Big)^2\mathbf{1}_{\{\z^*_0\leq a\}}\Bigg).$$
In order to prove (ii), we thus need to get the bound $J(a)\leq C\sqrt{a}$ when $a\geq 1$. 

To this end, we apply Theorem \ref{re-root-theo} with 
$$G(\omega)=\mathbf{1}_{\{\omega(0)<1\}}\,\mathbf{1}_{\{\Theta(\omega)\leq a\}}\,\int_0^\sigma\dd u\,\mathbf{1}_{\{\wh W_u(\omega)<1\}},$$
where the function $\Theta$ is defined on $\S$ by
$$\Theta(\omega)=\liminf_{\ve \to 0} \frac{1}{\ve^2}\int_0^{\sigma(\omega)} \dd s\,\mathbf{1}_{\{\hat\omega_s<\ve\}},$$
in such a way that we have $\z^*_0=\Theta(\omega^{[r]})$ for every $r\in(0,\sigma)$, $\N^*_0$ a.e., and $\z_0=\Theta(\tr_0(\omega))$, $\N_b$ a.e.
for every $b>0$
(use \eqref{approxz*0} and \eqref{formu-exit}).
It follows that
$$J(a)=2\int_0^1 \dd b \,\N_b\Big(\z_0\,\mathbf{1}_{\{\z_0\leq a\}}\int_0^\sigma\dd u\,\mathbf{1}_{\{\wh W_u<1,\zeta_u\leq\tau_0(W_u)\}}\Big),$$
where we recall the notation $\tau_0(\w)=\inf\{t\in(0,\zeta_{(\w)}]:\w(t)=0\}$, for $\w\in\mathcal{W}$. 

Let us fix $b\in(0,1)$ and set
$$K(a,b)=\N_b\Big(\z_0\,\mathbf{1}_{\{\z_0\leq a\}}\int_0^\sigma\dd u\,\mathbf{1}_{\{\wh W_u<1,\zeta_u\leq\tau_0(W_u)\}}\Big).$$
Then $K(a,b)\leq e\,\wt K(a,b)$, where
$$\wt K(a,b)=\N_b\Big(\z_0\,e^{-\z_0/a}\int_0^\sigma\dd u\,\mathbf{1}_{\{\wh W_u<1,\zeta_u\leq\tau_0(W_u)\}}\Big).$$
Let $(L^0_s)_{0\leq s\leq \sigma}$ denote the exit local time from $(0,\infty)$
as defined in formula \eqref{formu-exit}, and recall that
$\z_0=L^0_\sigma$. Then,
\begin{align*}
\wt K(a,b)=&\;\N_b\Big(\int_0^\sigma\dd u\,\mathbf{1}_{\{\wh W_u<1,\zeta_u\leq\tau_0(W_u)\}}L^0_u\,e^{-L^0_\sigma/a}\Big)\\
&+ \N_b\Big(\int_0^\sigma\dd u\,\mathbf{1}_{\{\wh W_u<1,\zeta_u\leq\tau_0(W_u)\}}(L^0_\sigma-L^0_u)\,e^{-L^0_\sigma/a}\Big),
\end{align*}
and the two terms in the right-hand side are equal, by a simple time-reversal argument. Let us consider the second term,
and bound $e^{-L^0_\sigma/a}$ by $e^{-(L^0_\sigma-L^0_u)/a}$. Using the Markov property under $\N_b$ (cf. formula \eqref{Markov-excu}), we
get
\begin{equation}
\label{bd-1}
\wt K(a,b)\leq 2\,\N_b\Bigg(\int_0^\sigma\dd u\,\mathbf{1}_{\{\wh W_u<1,\zeta_u\leq\tau_0(W_u)\}}\E^*_{W_u}\Big(L^0_\sigma\,e^{-L^0_\sigma/a}\Big)\Bigg),
\end{equation}
where we note that the definition of the exit local time also makes sense under $\P^*_\w$ for any $\w\in\W_b$, see
\cite[Section V.1]{Zurich}. 

Let $\w\in \W_b$ such that $\zeta_{(\w)}\leq \tau_0(\w)$. For every $\lambda>0$, we compute
$$\E^*_\w[L^0_\sigma e^{-\lambda L^0_\sigma}]=-\frac{\dd}{\dd\lambda} \E^*_\w[e^{-\lambda L^0_\sigma}].$$
We use Lemma V.5 of \cite{Zurich}, which says that the evolution of the Brownian snake 
under $\P^*_\w$ is described by a Poisson measure $\mathcal{P}$ on $[0,\zeta_{(\w)}]\times \S$
with intensity $2\,\dd t\,\N_{\w(t)}(\dd \omega)$, in such a  way that
$$ \E^*_\w[e^{-\lambda L^0_\sigma}]=\E^*_{\w}\Big[ \exp\Big(-\lambda\int \mathcal{P}(\dd t\dd \omega)\,\z_0(\omega)\Big)\Big]
= \exp\Big(-2\int_0^{\zeta_{(\w)}} \dd t\,\N_{\w(t)}(1-e^{-\lambda \z_0})\Big).$$
Using also \eqref{Laplace-exit}, we get
$$
\E^*_\w[L^0_\sigma e^{-\lambda L^0_\sigma}]
=-\frac{\dd}{\dd\lambda}\Bigg(\exp\Bigg(-3\int_0^{\zeta_{(\w)}} \dd t\,\Big(\w(t)+\sqrt{3/2\lambda}\Big)^{-2}\Bigg)\Bigg).
$$
It follows that
\begin{align*}
&\E^*_\w[L^0_\sigma e^{-\lambda L^0_\sigma}]\\
& =3\sqrt{3/2}\,\lambda^{-3/2} \Bigg(\int_0^{\zeta_{(\w)}} \frac{\dd t}{(\w(t)+\sqrt{3/2\lambda})^3}\Bigg)
\exp\Bigg(-3\int_0^{\zeta_{(\w)}} \dd t\,\Big(\w(t)+\sqrt{3/2\lambda}\Big)^{-2}\Bigg)\Bigg).
\end{align*}
We take $\lambda=1/a$ and substitute the identity of the last display in \eqref{bd-1}. From the first moment formula \eqref{first-mom},
we get
\begin{align*}
\wt K(a,b)\leq 6\sqrt{\frac{3}{2}}\,a^{3/2} \,\int_0^\infty \dd t \,&\mathbf{E}_b\Bigg[\mathbf{1}_{\{t\leq\kappa_0,B_t<1\}}\Bigg(\int_0^{t} \frac{\dd s}{(B_s+\sqrt{3a/2})^3}\Bigg)\\
&\qquad\times \exp\Bigg(-3\int_0^{t} \dd u\,\Big(B_u+\sqrt{3a/2}\Big)^{-2}\Bigg)\Bigg],
\end{align*}
where $(B_t)_{t\geq 0}$ is a linear Brownian motion that starts from $b$ under the probability measure $\mathbf{P}_b$, and, for
every $x\in\R$, $\kappa_x=\inf\{t\geq 0:B_t=x\}$. To simplify
notation, let us set $C_1:=6\sqrt{3/2}$, and write $a'=\sqrt{3a/2}$. It follows that
\begin{align*}
&\wt K(a,b)\\
&\quad\leq C_1\,a^{3/2} \,\mathbf{E}_b\Bigg[\int_0^\infty \frac{\dd s}{(B_s+a')^3}\int_s^\infty \dd t \,\mathbf{1}_{\{t\leq\kappa_0,B_t<1\}}
\exp\Bigg(-3\int_0^{t} \dd u\,(B_u+a')^{-2}\Bigg)\Bigg]\\
&\quad=C_1\,a^{3/2} \,\mathbf{E}_{b+a'}\Bigg[\int_0^\infty \frac{\dd s}{(B_s)^3}\int_s^\infty \dd t \,\mathbf{1}_{\{t\leq\kappa_{a'},B_t<a'+1\}}
\exp\Bigg(-3\int_0^{t} \dd u\,(B_u)^{-2}\Bigg)\Bigg]\\
&\quad\leq C_1\,a^{3/2} \int_0^\infty \dd s\, \mathbf{E}_{b+a'}\Bigg[(B_s)^{-3}\exp\Bigg(-3\int_0^{s} \dd u\,(B_u)^{-2}\Bigg)\int_s^\infty \dd t\,\mathbf{1}_{\{t\leq\kappa_{a'},B_t<a'+1\}}\Bigg]\\
&\quad=C_1\,a^{3/2} \int_0^\infty \dd s\, \mathbf{E}_{b+a'}\Bigg[(B_s)^{-3}\exp\Bigg(-3\int_0^{s} \dd u\,(B_u)^{-2}\Bigg)\,\mathbf{1}_{\{s<\kappa_{a'}\}}\,\\
&\hspace{8cm}\times\mathbf{E}_{B_s}\Big[\int_0^\infty \dd t\,\mathbf{1}_{\{t\leq\kappa_{a'},B_t<a'+1\}}\Big]\Bigg]
\end{align*}
where we have applied the Markov property at time $s$. 
We then observe that, for 
every $x>a'$,
\begin{equation}
\label{bd-intermed}
\mathbf{E}_{x}\Big[\int_0^\infty \dd t\,\mathbf{1}_{\{t\leq\kappa_{a'},B_t<a'+1\}}\Big]\leq 1.
\end{equation}
Indeed, the strong Markov property shows that it is enough to prove this when $x=a'+1$, and then
$$\mathbf{E}_{a'+1}\Big[\int_0^\infty \dd t\,\mathbf{1}_{\{t\leq\kappa_{a'},B_t<a'+1\}}\Big]
= \mathbf{E}_{1}\Big[\int_0^{\kappa_0} \dd t\,\mathbf{1}_{\{B_t<1\}}\Big]=\E\Big[\int_0^1\dd u\,(R^{(2)}_u)^2\Big]=1,$$
where $R^{(2)}$ denotes a two dimensional Bessel process started at $0$, and the second equality follows from a classical Ray-Knight theorem.

Thanks to \eqref{bd-intermed}, we arrive at
\begin{equation}
\label{bd-intermed2}
\wt K(a,b)\leq C_1\, a^{3/2}\,\int_0^\infty \dd s\, \mathbf{E}_{b+a'}\Bigg[(B_s)^{-3}\exp\Bigg(-3\int_0^{s} \dd u\,(B_u)^{-2}\Bigg)\,\mathbf{1}_{\{s<\kappa_{a'}\}}\Bigg].
\end{equation}
At this stage, we use the absolute continuity relations between Brownian motion and Bessel processes
(see e.g. Section 2 of \cite{Yor} or \cite[Proposition 2.6]{LGW}) to get, for every $s>0$,
\begin{align*}
&\mathbf{E}_{b+a'}\Bigg[(B_s)^{-3}\exp\Bigg(-3\int_0^{s} \dd u\,(B_u)^{-2}\Bigg)\,\mathbf{1}_{\{s<\kappa_{a'}\}}\Bigg]\\
&\qquad=(b+a')^3\,\mathbf{E}_{b+a'}\Big[(R_s)^{-6}\,\mathbf{1}_{\{R_u>a', \forall u\in[0,s]\}}\Big],
\end{align*}
where $(R_t)_{t\geq 0}$ denotes a Bessel process of dimension $7$ that starts at $x$ under the probability measure $\mathbf{P}_x$. 
Recalling that $b\in(0,1)$ and $a\geq 1$, the right-hand side of the last display is bounded above by $C'\,\mathbf{E}_{b+a'}[(R_s)^{-3}\,\mathbf{1}_{\{R_u>a', \forall u\in[0,s]\}}]$, for
some constant $C'$ independent of $b$ and $a$. 
From \eqref{bd-intermed2},
we then get
\begin{equation}
\label{bdinter3}
\wt K(a,b)\leq C'C_1\, a^{3/2}\,\int_0^\infty \dd s\, \mathbf{E}_{b+a'}\Big[(R_s)^{-3}\,\mathbf{1}_{\{R_u>a', \forall u\in[0,s]\}}\Big].
\end{equation}
Set $\kappa^R_{a'}=\inf\{t\geq 0:R_t\leq a'\}\leq +\infty$. Then,
$$\mathbf{E}_{b+a'}\!\Big[(R_s)^{-3}\,\mathbf{1}_{\{R_u>a', \forall u\in[0,s]\}}\Big]\!
=\mathbf{E}_{b+a'}\!\Big[(R_s)^{-3}\,\mathbf{1}_{\{R_s>a'\}}\Big] \!- \mathbf{E}_{b+a'}\!\Big[(R_s)^{-3}\,\mathbf{1}_{\{R_s>a',\,\kappa^R_{a'}\leq s\}}\Big].$$
On one hand,
\begin{align*}
\int_0^\infty \dd s\,\mathbf{E}_{b+a'}\Big[(R_s)^{-3}\,\mathbf{1}_{\{R_s>a',\,\kappa^R_{a'}\leq s\}}\Big]
&=\mathbf{E}_{b+a'}\Big[\int_{\kappa^R_{a'}}^\infty \dd s\,(R_s)^{-3}\,\mathbf{1}_{\{R_s>a'\}}\Big]\\
&=\mathbf{P}_{b+a'}(\kappa^R_{a'}<\infty) \times \mathbf{E}_{a'}\Big[\int_0^\infty \dd s\,(R_s)^{-3}\mathbf{1}_{\{R_s>a'\}}\Big]\\
&=\Big(\frac{a'}{b+a'}\Big)^5 \times \mathbf{E}_{0}\Big[\int_0^\infty \dd s\,(R_s)^{-3}\mathbf{1}_{\{R_s>a'\}}\Big]
\end{align*}
and, on the other hand,
$$\int_0^\infty \dd s\,\mathbf{E}_{b+a'}\Big[(R_s)^{-3}\,\mathbf{1}_{\{R_s>a'\}}\Big]\leq \mathbf{E}_{0}\Big[\int_0^\infty \dd s\,(R_s)^{-3}\mathbf{1}_{\{R_s>a'\}}\Big].$$
By combining the last two displays, we get
$$\int_0^\infty \dd s\, \mathbf{E}_{b+a'}\Big[(R_s)^{-3}\,\mathbf{1}_{\{R_u>a', \forall u\in[0,s]\}}\Big]
\leq \Big(1- \Big(\frac{a'}{b+a'}\Big)^5\Big)\, \mathbf{E}_{0}\Big[\int_0^\infty \dd s\,(R_s)^{-3}\mathbf{1}_{\{R_s>a'\}}\Big],$$
and we have
$$\mathbf{E}_{0}\Big[\int_0^\infty \dd s\,(R_s)^{-3}\,\mathbf{1}_{\{R_s>a'\}}\Big]
=\int_{\R^7} \frac{\dd z}{|z|^3}\,G(z)\,\mathbf{1}_{\{|z|>a'\}}\leq \frac{C''}{a'},
$$
where $G(z)=c\,|z|^{-5}$ denotes the Green function of Brownian motion in $\R^7$, and $C''$ is a constant. Finally, 
since 
$$1- \Big(\frac{a'}{b+a'}\Big)^5\leq 5\times\frac{b}{b+a'}\leq \frac{5}{a'},$$
the bound \eqref{bdinter3} gives
$$\wt K(a,b)\leq 5C''C'C_1\,a^{3/2}\,(a')^{-2}=(10/3)C''C'C_1\sqrt{a},$$
since $(a')^2= 3a/2$.  We conclude that 
$J(a)\leq C\sqrt{a}$, with $C=(20/3)eC''C'C_1$. This completes the proof.
\endproof

\section{Proof of Theorem \ref{mainT}}

Let us write $\mathcal{M}(\R)$ for the space of all finite measures on $\R$, which is equipped with the topology of
weak convergence and the associated Borel $\sigma$-field. We define a transition kernel from
$(0,\infty)\times \R$ into $\mathcal{M}(\R)$ as follows. For $(t,y)\in (0,\infty)\times\R$, we use the notation $U^{\mathrm{br},t,y}$
for the bridge of duration $t$ from $0$ to $y$ associated with the stable L\'evy process with no negative jumps
and Laplace exponent $\frac{1}{2}\psi(\lambda)=\sqrt{2/3}\,\lambda^{3/2}$ (we again refer to \cite{FPY}
for the construction of this bridge). Let $\eta_k$, $k\in\N$,
be the sequence of jumps of $U^{\mathrm{br},t,y}$ ranked in nonincreasing order. We define $Q((t,y),\dd\mu)$
as the probability measure on $\mathcal{M}(\R)$ obtained as the distribution of
$$\sum_{k\in\N} \y(\omega_k)$$
where, conditionally on $U^{\mathrm{br},t,y}$, the random snake trajectories $\omega_k$
are independent, and, for every $k$, $\omega_k$ is distributed according to $\N^{*,\eta_k}_0$. This definition 
makes sense because, using formula \eqref{N*sigma},
$$\E\Big[\sum_{k\in\N} \langle\y(\omega_k),1\rangle\,\Big|\,U^{\mathrm{br},t,y} \Big]= \E\Big[\sum_{k\in\N}\sigma(\omega_k)
\,\Big|\,U^{\mathrm{br},t,y} \Big]=\sum_{k\in\N} (\eta_k)^2 <\infty, \quad\hbox{a.s.}$$
As usual, if $F$ is a nonnegative measurable function on $\mathcal{M}(\R)$, $QF$ stands for the function
on $(0,\infty)\times \R$ defined by
$$QF(t,y)= \int Q((t,y),\dd\mu)\,F(\mu).$$
We extend this definition by setting $QF(0,y)=F(0)$ for every $y\in\R$.

In order to prove Theorem \ref{mainT}, we will now argue under the measure $\N_0$. 
Recall the definition \eqref{occu-measure} of the random measure $\y$ and, for every $h\in\R$, 
let $\y^h_-$, resp. $\y^h_+$, denote the restriction of $\y$
to $(-\infty,h)$, resp. to $(h,\infty)$. We also let $\wt\y^h_+$, resp. $\wt\y^h_-$, be the pushforward of $\y^h_+$,
resp. of $\y^h_-$,
under the mapping $x\mapsto x-h$.
The key to the proof of Theorem \ref{mainT} is the following proposition.

\begin{proposition}
\label{keymainT}
Let $h\geq 0$. Let 
$F_1$
 and $F_2$ be two nonnegative measurable functions on $\mathcal{M}(\R)$. Then,
 $$\N_0\Big(F_1(\y^h_-)\,F_2(\wt\y^h_+)\Big)= \N_0\Big(F_1(\y^h_-)\,QF_2(\ell^h,\frac{1}{2}\dot\ell^h)\Big).$$
 \end{proposition}

Both assertions of Theorem \ref{mainT} follow from Proposition \ref{keymainT}. Just note that $(\ell^{h+x})_{x>0}$ is the
(continuous) density of the measure $\wt\y^h_+$, so that Proposition \ref{keymainT} immediately shows that 
the process $(\ell^{h+x},\dot\ell^{h+x})_{x\geq 0}$ is independent of $(\ell^{h+x},\dot\ell^{h+x})_{x\leq 0}$ conditionally
on $(\ell^h,\dot\ell^h)$, and moreover its conditional distribution does not depend on $h$. 
The second assertion of Theorem \ref{mainT} follows from the case $h=0$ of Proposition \ref{keymainT}.

\proof[Proof of Proposition \ref{keymainT}]
We will use the fact that the local time $\ell^h$
can be expressed in terms of the exit measure process $(\mathcal{X}^h_r)_{r\geq0}$
via the formula
\begin{equation}
\label{local3}
\ell^h=\int_0^\infty \dd r \,\XX^h_r.
\end{equation}
See \cite[Proposition 26]{LGR} when $h>0$, and \cite[Proposition 3.1]{LGR2} when $h=0$.

In most of the proof, we deal with the case $h>0$, and, only at the end, we explain how the desired result can be extended to 
the case $h=0$. So, we fix $h>0$, and we first note that the event $\ell^h=0$ occurs if and only if $\z_h=0$ (by \eqref{local3}
and the fact that $\XX^h_0=\z_h$). On
the event $\{\z_h=0\}$, we have $\wt\y^h_+=0$ and $QF_2(\ell^h,\dot\ell^h)=F_2(0)$. Thanks to this observation,
it is enough to prove the formula of the proposition with
$\N_0$ replaced by the conditional
probability measure $\N^{(h)}_0(\dd\omega):=\N_0(\dd\omega\mid \z_h>0)=\N_0(\dd\omega\mid W^*>h)$. 
From now on until the end of the proof when we consider the case $h=0$, we argue under the probability measure $\N^{(h)}_0$.
Recall from Section \ref{sec:exitlocal} that (under $\N^{(h)}_0$) the process $\XX^h$
is independent of $\tr_h(\omega)$ conditionally on $\z_h$.

We rely on the excursion theory presented in 
Section \ref{sec:excu} above, and we consider the excursions above and below level $h$, which are denoted 
by $\omega^{(u_i)}$, $i\in\N$, in Section \ref{sec:excu}. To simplify notation, we write $\omega^{(i)}$ instead of $\omega^{(u_i)}$
in this proof.
Recall that each excursion $\omega^{(i)}$ corresponds to
a jump time $r_i$ of the exit measure process $\XX^h$. We write $\delta_i:=\Delta\XX^h_{r_i}$ for
the corresponding jump. The conditional distribution of the collection $(\omega^{(i)})_{i\in\N}$ knowing the process $\XX^h$
(and $\tr_h(\omega)$)
is given by Theorem \ref{theo-excursion}.

For every $i\in \N$, let $\eta_i=1$ if $\omega^{(i)}$ is an excursion above $h$ and $\eta_i=-1$ otherwise. Notice that,
conditionally on the exit measure process $\XX^h$ (and on $\tr_h(\omega)$), the random variables $\eta_i$, $i\in \N$, 
are independent and uniformly distributed on $\{-1,+1\}$. Set $I_+:=\{i\in \N:\eta_i=+1\}$
and $I_-:=\{i\in \N:\eta_i=-1\}$.  As we already observed in  Section \ref{sec:excu}, if $\mathcal{C}$ is a connected
component of $\{v\in\t(\omega):V_v\not =h\}$ (other than the component containing the root),
and $u$ is the associated excursion debut,
the restriction of the volume measure $\mathrm{Vol}_{(\omega)}$ to $\mathcal{C}$ is identified with the volume measure $\mathrm{Vol}_{(\omega^{(u)})}$
(with the notation of Section \ref{sec:excu}). Similarly, it is not
hard to check that the restriction of $\mathrm{Vol}_{(\omega)}$
to the component containing the root is identified with the volume measure $\mathrm{Vol}_{(\mathrm{tr}_h(\omega))}$.
It follows from these remarks that
\begin{equation}
\label{local1}
\y^h_-= \y^h_-(\tr_h(\omega)) + \sum_{i\in I_-} \y^h_-(\omega^{(i)})
\end{equation}
and
\begin{equation}
\label{local2}
\y^h_+= \sum_{i\in I_+} \y^h_+(\omega^{(i)}).
\end{equation}

By the classical Lamperti transformation \cite{Lam}, we can write 
$(\XX^h_r)_{r\geq 0}$ as a time change of a 
L\'evy process stopped upon hitting $0$.
More precisely, we have for every $r\geq 0$,
$$\XX^h_r=\mathcal{U}_{\int_0^r \dd t \,\XX^h_t},$$
where $(\mathcal{U}_t)_{0\leq t\leq T_0}$ is 
a stable L\'evy process with no negative jumps and Laplace exponent $\psi$,
which is started at $\mathcal{U}_0=\z_h$
and stopped at its first hitting time of $0$.
%
Note that we have in particular
$$\int_0^\infty \dd r\,\XX^h_r = \inf\{t\geq 0:\mathcal{U}_t=0\}=T_0.$$
Recalling \eqref{local3}, we have thus $\ell^h=T_0$. 
We observe that $(\mathcal{U}_t)_{0\leq t\leq T_0}$ has the same jumps as $\XX^h$. Hence, for every $i\in \N$,
$\delta_i$ is the jump of $\mathcal{U}$ occurring at a certain time $s_i\in[0,T_0]$. 

The values of the process $\mathcal{U}$
 are determined by $(\XX^h_r)_{r\geq 0}$ only up to time $T_0$. With a small 
 abuse of notation, we can assume that the L\'evy process $(\mathcal{U}_t)_{t\geq 0}$
 is defined at all times under the underlying probability measure $\N^{(h)}_0(\dd\omega)$ (and is independent of 
 $\tr_h(\omega)$ conditionally on $\z_h$). 
 Let $(t_j)_{j\in \N}$ be the jump times of $\mathcal{U}$ (listed according to some 
 measurable enumeration) and, for every $j\in \N$, let $\gamma_j=\Delta \mathcal{U}_{t_j}$
 be the corresponding jump. Notice that, if $j\in \N$ is such that 
 $t_j\leq T_0$, there exists a unique $i\in \N$ such that $t_j=s_i$ and $\gamma_j=\delta_i$. 
  
 We may also assume that we have assigned a random snake trajectory
 $\ov\omega^j$ to each jump time $t_j$ of $\mathcal{U}$, in such a way that, if $t_j\leq T_0$, we have $\ov\omega^j=\omega^{(i)}$ where $i\in \N$ 
 is the unique index such that $t_j=s_i$, and,
 conditionally on $\mathcal{U}$,
 the random variables $\ov \omega^j$, $j\in \N$, are independent 
 and 
 the conditional distribution of $\ov\omega^j$ is
 $$\frac{1}{2}\N^{*,\gamma_j}_h + \frac{1}{2}\check\N^{*,\gamma_j}_h.$$
If $j\in \N$, we set $\ve_j=+1$ if 
 $\ov\omega^j$ is an excursion above $h$ and $\ve_j=-1$ otherwise. 
 We note that the ``labels'' $\ve_j$, $j\in\N$, are independent and 
 uniformly distributed over $\{-1,1\}$ (and are also independent of
 the process $\mathcal{U}$).
Set $J_+:=\{j\in \N:\ve_j=+1\}$ and $J_-:=\{j\in \N:\ve_j=-1\}$.

Let $\mathcal{U}'$ be the (centered) L\'evy process that is obtained from $\mathcal{U}$ by
``keeping only'' the jumps with label $-1$. More precisely, noting that the L\'evy measure of $\mathcal{U}$
is the measure $\bn(\dd z)$ defined in \eqref{Levy-mea}, we have for every
$t\geq 0$,
$$\mathcal{U}'_t=\z_h + \lim_{\alpha\downarrow 0} \Bigg(\sum_{j\in J_-,t_j\leq t,\gamma_j>\alpha} \gamma_j - \frac{t}{2} \int_\alpha^\infty x\,\bn(\dd x)\Bigg).$$
We also define $\mathcal{U}''_t=\mathcal{U}_t- \mathcal{U}'_t$, so that
\begin{equation}
\label{Levy-posi}
\mathcal{U}''_t=\lim_{\alpha\downarrow 0} \Bigg(
\sum_{j\in J_+,t_j\leq t,\gamma_j>\alpha} \gamma_j - \frac{t}{2} \int_\alpha^\infty x\,\bn(\dd x)\Bigg).
\end{equation}
Observe that $\mathcal{U}'$ and $\mathcal{U}''$ are two independent (centered) L\'evy processes with Laplace exponent $\frac{1}{2}\psi$,
such that $\mathcal{U}'_0=\z_h$ and $\mathcal{U}''_0=0$,  and also note that $\mathcal{U}'_{T_0}+\mathcal{U}''_{T_0}=\mathcal{U}_{T_0}=0$. 

At this point, it will be convenient to condition on the value of $\z_h$, and,
for every $z>0$, we introduce the conditional probability measure $\P^{(z)}:=\N^{(h)}_0(\cdot \mid\z_h=z)$,
in such a way that $\mathcal{U}'_0=\mathcal{U}_0=z$, $\P^{(z)}$ a.s. 
Then let $g$ be a nonnegative measurable
function on $\R_+$, and let $G_1$ and $G_2$ be two nonnegative measurable
functions on $\R_+\times \D(\R_+,\R)$, such that, for every $t\geq 0$, the mapping 
$\w\mapsto G_1(t,\w)$ is a (measurable) function of the collection 
of jumps of $(\w(s),0\leq s\leq t)$, and 
similarly for $G_2$. If $\pi_z(\dd t)$ denotes the law of $T_0$ under $\P^{(z)}$,
we have then 
\begin{equation}
\label{key00}
\E^{(z)}[g(T_0)G_1(T_0,\mathcal{U}')G_2(T_0,\mathcal{U}'')]= \int \pi_z(\dd t)\,g(t)\,\E^{(z)}[G_1(t,\mathcal{U}')G_2(t,\mathcal{U}'')\mid T_0=t].
\end{equation}
By Lemma \ref{Kemp} and the subsequent remarks, we know that the conditional distribution of the collection of 
jumps of $(\mathcal{U}_t)_{0\leq t\leq T_0}$ knowing that $T_0=t$ is 
the distribution of the collection of jumps of the bridge of duration $t$ from $z$ to $0$ associated with the L\'evy process $U$.
Recalling that the signs 
$\ve_j$ are assigned independently knowing $\mathcal{U}$, it follows that, for every $t\geq 0$,
\begin{equation}
\label{key0}
\E^{(z)}[G_1(t,\mathcal{U}')G_2(t,\mathcal{U}'')\mid T_0=t]=\E^{(z)}[G_1(t,\mathcal{U}')G_2(t,\mathcal{U}'')\mid \mathcal{U}_t=0].
\end{equation}
Now note that $\mathcal{U}_t=0$ is equivalent to $\mathcal{U}''_t=-\mathcal{U}'_t$. Using the independence of $\mathcal{U}'$ and $\mathcal{U}''$, we 
can verify that
\begin{equation}
\label{key1}
\E^{(z)}[G_1(t,\mathcal{U}')G_2(t,\mathcal{U}'')\mid \mathcal{U}_t=0]=\E^{(z)}[G_1(t,\mathcal{U}')\Phi(t,-\mathcal{U}'_t)\mid \mathcal{U}_t=0],
\end{equation}
where we use the notation $\Phi(t,a):=\E^{(z)}[G_2(t,\mathcal{U}'')\mid \mathcal{U}''_t=a]$ for every $a\in \R$
(this function does not depend on $z$).
The identity \eqref{key1} may be derived from elementary manipulations. Alternatively, we 
may proceed as follows. We set $\wt{\mathcal{U}}_s=\mathcal{U}'_s$ if $s\in[0,t]$ and $\wt{\mathcal{U}}_s= \mathcal{U}'_t + \mathcal{U}''_{s-t}$ if $s\in [t,2t]$,
so that, under $\P^{(z)}(\cdot\mid \mathcal{U}_t=0)$, $(\wt{\mathcal{U}}_s)_{s\in[0,2t]}$ is a L\'evy process conditioned on $\wt{\mathcal{U}}_{2t}=0$.
Then \eqref{key1} is nothing but the usual Markov property at time $t$ for the L\'evy process bridge. 

Thanks to \eqref{key0} and \eqref{key1},
we get 
$$\E^{(z)}[G_1(t,\mathcal{U}')G_2(t,\mathcal{U}'')\mid T_0=t]=\E^{(z)}[G_1(t,\mathcal{U}')\Phi(t,-\mathcal{U}'_t)\mid T_0=t],$$
where we also use the fact that $\mathcal{U}'_t$ is a measurable function of the jumps of $\mathcal{U}'$
over $[0,t]$. Recalling \eqref{key00}, we finally get that
\begin{align}
\label{key2}
\E^{(z)}[g(T_0)G_1(T_0,\mathcal{U}')G_2(T_0,\mathcal{U}'')]&= \int \pi_z(\dd t)\,g(t)\,\E^{(z)}[G_1(t,\mathcal{U}')\Phi(t,-\mathcal{U}'_t)\mid T_0=t]\\\
&=\E^{(z)}[g(T_0)G_1(T_0,\mathcal{U}')\Phi(T_0,-\mathcal{U}'_{T_0})].\nonumber
\end{align}

We now use Theorem \ref{theo-excursion}.
Under $\N_0^{(h)}$, conditionally on $(\mathcal{X}^h_r)_{r\geq 0}$,
the excursions $ \omega^{(i)}$, $i\in \N$, are independent, and independent
of $\tr_h(\omega)$, and the 
conditional distribution of $\omega^{(i)}$ is 
$\frac{1}{2}(\N^{*,\delta_i}_0+\check\N^{*,\delta_i}_0)$.
It easily follows that, conditionally on $(\mathcal{X}^h_r)_{r\geq 0}$ and on $(\eta_i)_{i\in\N}$,
the excursions $ \omega^{(i)}$, $i\in \N$, are independent (and independent
of $\tr_h(\omega)$), and the 
conditional distribution of $\omega^{(i)}$ is $\N^{*,\delta_i}_h$ if $\eta_i=1$ (equivalently $i\in I_+$)
and $\check\N^{*,\delta_i}_h$ if $\eta_i=-1$ (equivalently $i\in I_-$).
Hence, if $H,A,B$ are nonnegative measurable functions on the space of snake trajectories, we get
\begin{align*}
&\N^{(h)}_0\Bigg(H(\tr_h(\omega))\,\exp\Big(-\sum_{i\in I_-} A(\omega^{(i)})\Big)\exp\Big(-\sum_{i\in I_+} B(\omega^{(i)})\Big)\Bigg)\\
&\qquad= \N^{(h)}_0\Bigg(H(\tr_h(\omega))\,\prod_{i\in I_-} \check\N^{*,\delta_i}_h(e^{-A})\,\prod_{i\in I_+} \N^{*,\delta_i}_h(e^{-B})\Bigg)\\
&\qquad= \N^{(h)}_0\Bigg(H(\tr_h(\omega))\,\E^{(\z_h)}\Bigg[\prod_{i\in I_-} \check\N^{*,\delta_i}_h(e^{-A})\,\prod_{i\in I_+} \N^{*,\delta_i}_h(e^{-B})\,\Bigg]\Bigg)
\end{align*}
The quantities
$$\prod_{i\in I_-} \check\N^{*,\delta_i}_h(e^{-A}),\quad \prod_{i\in I_+} \N^{*,\delta_i}_h(e^{-B})$$
are functions of the jumps of $\mathcal{U}'$ and $\mathcal{U}''$, respectively, over the time interval $[0,T_0]$. Hence, we can use 
\eqref{key2} to get, for every $z>0$,
$$\E^{(z)}\Bigg[\prod_{i\in I_-} \check\N^{*,\delta_i}_h(e^{-A})\,\prod_{i\in I_+} \N^{*,\delta_i}_h(e^{-B})\,\Bigg]
= \E^{(z)}\Bigg[\prod_{i\in I_-} \check\N^{*,\delta_i}_h(e^{-A})\,\Phi_B(T_0,-\mathcal{U}'_{T_0})\Bigg],$$
where $\Phi_B(t,y)$ is the expected value of the quantity
$$\prod_{k\in \N} \N^{*,a_k}_h(e^{-B})$$
where the numbers $a_k$, $k\in \N$ are the jumps of the bridge of duration $t$ from $0$ to $y$,
for a L\'evy process with no negative jumps
and Laplace exponent $\frac{1}{2}\psi$. 
We finally conclude that
\begin{align}
\label{key3}
&\N^{(h)}_0\Bigg(H(\tr_h(\omega))\,\exp\Big(-\sum_{i\in I_-} A(\omega^{(i)})\Big)\exp\Big(-\sum_{i\in I_+} B(\omega^{(i)})\Big)\Bigg)\\
&\qquad=\N^{(h)}_0\Bigg(H(\tr_h(\omega))\,\prod_{i\in I_-} \check\N^{*,\delta_i}_h(e^{-A})\,\Phi_B(T_0,-\mathcal{U}'_{T_0})\Bigg)\nonumber\\
&\qquad=\N^{(h)}_0\Bigg(H(\tr_h(\omega))\,\exp\Big(-\sum_{i\in I_-} A(\omega^{(i)})\Big)\,\Phi_B(T_0,-\mathcal{U}'_{T_0})\Bigg).\nonumber
\end{align}

\begin{lemma}
\label{ident-derivative} $\N^{(h)}_0$ a.s., 
we have $\mathcal{U}'_{T_0}=-\mathcal{U}''_{T_0}=-\frac{1}{2}\dot\ell^h$, where $\dot\ell^h$ denotes the derivative at $h$
of the function $x\mapsto \ell^x$.
\end{lemma}

Let us postpone the proof of Lemma \ref{ident-derivative}. 
Since we already know that $T_0=\ell^h$, we have $\Phi_B(T_0,-\mathcal{U}'_{T_0})=\Phi_B(\ell^h, \frac{1}{2}\dot\ell^h)$
in formula \eqref{key3}. Next let $f_1$ and $f_2$ be two bounded measurable functions on $\R$, and
consider the functions $F_1$ and $F_2$ defined on $\mathcal{M}(\R)$ by $F_i(\mu)=\exp-\langle \mu,f_i\rangle$, for $i=1,2$. 
Recalling \eqref{local1} and \eqref{local2}, we see that an appropriate choice of the 
functions $H,A,B$ in \eqref{key3}  gives
$$\N^{(h)}_0\Big(F_1(\y^h_-)\,F_2(\wt\y^h_+)\Big)= \N^{(h)}_0\Big(F_1(\y^h_-)\,\Phi_{(f_2)}(\ell^h, \frac{1}{2}\dot\ell^h)\Big),$$
where $\Phi_{(f_2)}(t,y)=QF_2(t,y)$, with the notation introduced before Proposition \ref{keymainT}. 
We have thus obtained the special case of the formula of Proposition \ref{keymainT} when $F_1$ and $F_2$ are as specified above, and 
a standard monotone class argument (see e.g. Lemma II.5.2 in \cite{Per}) gives the general case. 
This completes the proof in the case $h>0$. 

Consider now the case $h=0$. It seems plausible that one could derive this case 
by passing to the limit $h\to 0$ in the formula obtained for $h>0$. However, a rigorous justification
of this passage to the limit leads to certain technical difficulties, and, for this reason, we will
use a different argument based on the re-rooting property of $\N_0$. For $\omega\in\S_0$ and 
$r\in[0,\sigma(\omega)]$, recall the notation 
$\wt\omega^{[r]}$ introduced before formula \eqref{re-root-for}, and note that we have $\N_0$ a.e.
$$\y_-^0(\wt\omega^{[r]})\!=\wt\y_-^{\hat\omega_r}(\omega),\;\y_+^0(\wt\omega^{[r]})\!=\wt\y_+^{\hat\omega_r}(\omega),\; \ell^0(\wt\omega^{[r]})\!=\ell^{\hat\omega_r}(\omega),
\;\dot\ell^0(\wt\omega^{[r]})\!=\dot\ell^{\hat\omega_r}(\omega),\;\wh{\wt\omega}^{[r]}_{\sigma-r}\!=-\wh\omega_r.$$
Let $F_1$ and $F_2$ be nonnegative measurable functions on $\mathcal{M}(\R)$. From formula \eqref{re-root-for} and the preceding display, we get 
\begin{align}
\label{techh0}
&\N_0\Bigg( \int_0^\sigma \dd r\,\mathbf{1}_{\{\hat\omega_r>0\}}\,F_1(\wt\y_-^{\hat\omega_r}(\omega))\,F_2(\wt\y_+^{\hat\omega_r}(\omega))\Bigg)\\
&\qquad=
\N_0\Bigg( \int_0^\sigma \dd r\,\mathbf{1}_{\{\hat\omega_r>0\}}\,F_1(\y^0_-(\wt\omega^{[r]}))\,F_2(\y^0_+(\wt\omega^{[r]}))\Bigg)\nonumber\\
&\qquad=\N_0\Bigg(\int_0^\sigma \dd r\,\mathbf{1}_{\{\hat\omega_{\sigma-r}<0\}}\,F_1(\y^0_-(\omega))\,F_2(\y^0_+(\omega))\Bigg)\nonumber\\
&\qquad=\N_0\Big(\langle \y^0_-,1\rangle\,F_1(\y^0_-)\,F_2(\y^0_+)\Big).\nonumber
\end{align}
On the other hand, the left-hand side of \eqref{techh0} is also equal to
\begin{align}
\label{techh1}
\N_0\Bigg(\int \y^0_+(\dd x)\,F_1(\wt\y^x_-)\,F_2(\wt\y^x_+)\Bigg)&=\N_0\Bigg(\int_0^\infty \dd x\,\ell^x\,F_1(\wt\y^x_-)\,F_2(\wt\y^x_+)\Bigg)\\
&=\int_0^\infty \dd x\,\N_0\Big(\ell^x\,F_1(\wt\y^x_-)\,F_2(\wt\y^x_+)\Big)\nonumber\\
&=\int_0^\infty \dd x\,\N_0\Big(\ell^x\,F_1(\wt\y^x_-)\,QF_2(\ell^x,\frac{1}{2}\dot\ell^x)\Big),\nonumber
\end{align}
where we use the case $h>0$ of Proposition \ref{keymainT} in the last equality. Finally, replacing the function $F_1(\mu)$
by  $(\langle \mu,1\rangle)^{-1}F_1(\mu)$, we deduce from \eqref{techh0} and \eqref{techh1} that
$$\N_0\Big(F_1(\y^0_-)\,F_2(\y^0_+)\Big)= \int_0^\infty \dd x\,\N_0\Big(\ell^x\,(\langle \wt\y^x_-,1\rangle)^{-1}\,F_1(\wt\y^x_-)\,QF_2(\ell^x,
\frac{1}{2}\dot\ell^x)\Big).$$
The right-hand side of the preceding display remains the same if we take $F_2=1$ and replace $F_1(\y^0_-)$ by $F_1(\y^0_-)QF_2(\ell^0,\frac{1}{2}\dot\ell^0)$:
Note that the pair $(\ell^0,\frac{1}{2}\dot\ell^0)$ is a measurable function of $\y^0_-$, such that the same function applied to 
the measure $\wt\y^x_-$ gives $(\ell^x,\frac{1}{2}\dot\ell^x)$. The case $h=0$ of Proposition \ref{keymainT} now follows. \endproof

\proof[Proof of Lemma \ref{ident-derivative}]
To simplify notation, we write (only in this proof) $\P$ for the probability measure $\N^{(h)}_0$
and $\E$ for the corresponding expectation.
We have already noted that  $\mathcal{U}'_{T_0}=-\mathcal{U}''_{T_0}$, and so we only need to verify that $\mathcal{U}''_{T_0}=\frac{1}{2}\dot\ell^h$.
We first observe that
\begin{align}
\label{conv-deri}
\frac{1}{\ve^2} \Bigg(\int_0^\sigma \dd s\,\mathbf{1}_{\{h<\wh W_s<h+\ve\}} - \ve\,\ell^h\Bigg)
=\frac{1}{\ve^2} \Bigg(\int_h^{h+\ve} \dd x\,\ell^x -\ve \ell^h\Bigg)
&=\int_0^1 \dd y \Big(\frac{\ell^{h+\ve y}-\ell^h}{\ve}\Big)\\
&\build{\la}_{\ve\to 0}^{}  \frac{1}{2}\,\dot\ell^h,\nonumber
\end{align}
$\P$ a.s. On the other hand, we have
\begin{align}
\label{time-below}
\int_0^\sigma \dd s\,\mathbf{1}_{\{h<\wh W_s<h+\ve\}} &= \sum_{i\in I_+} \int_0^{\sigma(\omega^{(i)})}\dd s\,\mathbf{1}_{\{h<\wh W_s(\omega^{(i)})<h+\ve\}}\\
&=\sum_{j\in J_+,t_j\leq T_0} \int_0^{\sigma(\ov\omega^j)}\dd s\,\mathbf{1}_{\{h<\wh W_s(\ov\omega^j)<h+\ve\}}.\nonumber
\end{align}
For every $j\in J_+$ and $\ve>0$, set
$$\gamma_j^\ve = \frac{1}{\ve ^2}\int_0^{\sigma(\ov\omega_j)}  \dd s\,\mathbf{1}_{\{h<\wh W_s(\ov\omega^j)<h+\ve\}}.$$
Recall that, conditionally on $\XX^h$ and on $\{j\in J_+\}$, $\ov\omega_j$ is distributed according 
to $\N^{*,\gamma_j}_h$. 
By  \eqref{approxz*0} and the remarks following the definition of $\N^{*,z}_0$, we have $\gamma_j^\ve\la \gamma_j=\Delta \mathcal{U}_{t_j}$ as $\ve\to 0$, for
every $j\in J_+$,  $\P$ a.s.

Let $\alpha>0$. Since the set $\{j\in J_+:\gamma_j\geq \alpha,t_j\leq T_0\}$ is finite, it follows that
\begin{equation}
\label{large-jump}
\sum_{j\in J_+,\gamma_j\geq \alpha,t_j\leq T_0} \gamma^\ve_j \build{\la}_{\ve\to 0}^{} \sum_{j\in J_+,\gamma_j\geq \alpha,t_j\leq T_0} \gamma_j,
\quad\P\hbox{ a.s.}
\end{equation}

%

%
For every $\ve>0$ and $0\leq u<v\leq \infty$, we set
$$\Gamma_\ve(u,v)= \ve^{-2} \int_u^v \bn(\dd z)\,\N^{*,z}_0\Bigg(\int_0^\sigma \dd s\,\mathbf{1}_{\{0<\wh W_s<\ve\}}\Bigg)=\ve^2 \int_u^v f(\frac{z}{\ve^2})\,\bn(\dd z),$$
with the notation of Lemma \ref{tech-lemma}.
We observe that
$$\Gamma_\ve(0,\infty)=\ve^{-2}\N^*_0\Bigg(\int_0^\sigma \dd s\,\mathbf{1}_{\{\wh W_s<\ve\}}\Bigg)=\frac{2}{\ve}$$
by \eqref{first-mo}. Moreover,
\begin{equation}
\label{large-jump2}
\Gamma_\ve(\alpha,\infty)
=\int_\alpha^\infty \ve^2\,f(\frac{z}{\ve^2})\,\bn(\dd z)\build{\la}_{\ve \to 0}^{} \int_\alpha^\infty z\,\bn(\dd z),
\end{equation}
 by dominated convergence (justified by Lemma \ref{tech-lemma} (i)). 

By construction (and standard properties of L\'evy processes), the point measure
\begin{equation}
\label{Poisson-mea}
\sum_{j\in J_+} \delta_{(t_j,\gamma_j,\ov\omega^j)}
\end{equation}
is Poisson with intensity $\dd t\,\frac{1}{2}\bn(\dd z)\,\N^{*,z}_h(\dd \omega)$. In particular, for $0\leq u<v\leq \infty$, we have
$$\E\Bigg[\sum_{j\in J_+,u\leq \gamma_j<v,t_j\leq t} \gamma^\ve_j \Bigg]= \frac{t}{2} \Gamma_\ve(u,v).$$
Then, using a classical formula for Poisson measures
(see formula (3.19) in \cite{King}), we have
\begin{align*}&\E\Bigg[\Bigg(\sum_{j\in J_+,\gamma_j< \alpha,t_j\leq t} \gamma^\ve_j -\E\Bigg[\sum_{j\in J_+,\gamma_j< \alpha,t_j\leq t} \gamma^\ve_j\Bigg]\Bigg)^2\Bigg]\\
&\qquad= \frac{t}{2} \int_0^\alpha \bn(\dd z)\,\N^{*,z}_h\Bigg( \Big(\ve^{-2}\int_0^\sigma \dd s\,\mathbf{1}_{\{\wh W_s\leq h+\ve\}}\Big)^2\Bigg)\leq C t\sqrt{\alpha},
\nonumber
\end{align*}
where the last bound holds by Lemma \ref{tech-lemma} (ii) provided that $\alpha\geq \ve^2$. Under the latter condition, we can apply Doob's inequality in $L^2$ to the martingale
$$
M^\ve_t:=\sum_{j\in J_+,\gamma_j< \alpha,t_j\leq t} \gamma^\ve_j -\E\Bigg[\sum_{j\in J_+,\gamma_j< \alpha,t_j\leq t} \gamma^\ve_j\Bigg]
= \sum_{j\in J_+,\gamma_j< \alpha,t_j\leq t} \gamma^\ve_j -\frac{t}{2}\Gamma_\ve(0,\alpha)$$
and we get, for every $K>0$ and $\ve\in(0,\sqrt{\alpha}\,]$,
\begin{equation}
\label{small-jump}
\E\Bigg[\sup_{t\in[0,K]} \Bigg(\sum_{j\in J_+,\gamma_j<\alpha,t_j\leq t} \gamma^\ve_j -\frac{t}{2}\,\Gamma_\ve(0,\alpha)\Bigg)^2\Bigg]
\leq 4 CK\sqrt{\alpha}.
\end{equation}

Let us fix $\beta>0$. We observe that the convergence in \eqref{Levy-posi} holds uniformly when $t$ varies in a compact
set, at least along a suitable sequence of values of $\alpha$ decreasing to $0$ (see e.g. the proof of Theorem 1 in Chapter 1 of \cite{Be}). So we can choose
$\alpha>0$ small enough so that
\begin{equation}
\label{final1}
\P\Bigg(\Bigg|\Bigg(\sum_{j\in J_+,\gamma_j\geq\alpha,t_j\leq T_0} \gamma_j - \frac{T_0}{2} \int_\alpha^\infty z\,\bn(\dd z)\Bigg) -\mathcal{U}''_{T_0}\Bigg|>\beta\Bigg)<\beta.
\end{equation}
By choosing $\alpha$ even smaller if necessary, we may also assume thanks to \eqref{small-jump}
that, for every $\ve\in(0,\sqrt{\alpha}\,]$,
\begin{equation}
\label{final2}
\P\Bigg(\Bigg|\sum_{j\in J_+,\gamma_j< \alpha,t_j\leq T_0} \gamma^\ve_j -\frac{T_0}{2}\,\Gamma_\ve(0,\alpha)\Bigg|
>\beta\Bigg)<\beta.
\end{equation}
Once we have fixed $\alpha$, we can use \eqref{large-jump} and \eqref{large-jump2} to get that, for every small enough $\ve >0$, we have 
\begin{equation}
\label{final3}
\P\Bigg(\Bigg| \Big(\sum_{j\in J_+,\gamma_j\geq \alpha,t_j\leq T_0} \gamma^\ve_j-\frac{T_0}{2}\,\Gamma_\ve(\alpha,\infty)\Big)
- \Big(\sum_{j\in J_+,\gamma_j\geq \alpha,t_j\leq T_0} \gamma_j - \frac{T_0}{2} \int_\alpha^\infty z\,\bn(\dd z)\Big)\Bigg| >\beta\Bigg) 
<\beta.
\end{equation}
By combining \eqref{final1}, \eqref{final2} and \eqref{final3}, and using $\Gamma_\ve(0,\infty)=\Gamma_\ve(0,\alpha)+\Gamma_\ve(\alpha,\infty)$, we obtain that,
for $\ve$ small,
$$\P\Bigg(\Bigg|\Bigg(\sum_{j\in J_+,t_j\leq T_0} \gamma^\ve_j- \frac{T_0}{2}\,\Gamma_\ve(0,\infty)\Bigg)- \mathcal{U}''_{T_0}\Bigg|>3\beta\Bigg)
<3\beta.$$

Since $\beta$ was arbitrary, we have proved that
$$\sum_{j\in J_+,t_j\leq T_0} \gamma^\ve_j- \frac{T_0}{2}\,\Gamma_\ve(0,\infty) \build{\la}_{\ve\to 0}^{} \mathcal{U}''_{T_0}$$
in probability. Now recall from \eqref{time-below} that 
$$\sum_{j\in J_+,t_j\leq T_0} \gamma^\ve_j=\ve^{-2}\int_0^\sigma \dd s\,\mathbf{1}_{\{h<\wh W_s<h+\ve\}}.$$ Since we 
have also $\Gamma_\ve(0,\infty)=2/\ve$ and $T_0=\ell^h$, we conclude that
$$\frac{1}{\ve^2}\int_0^\sigma \dd s\,\mathbf{1}_{\{h<\wh W_s<h+\ve\}} - \frac{\ell^h}{\ve} \build{\la}_{\ve\to 0}^{}\mathcal{U}''_{T_0}$$
in probability. 
Comparing with \eqref{conv-deri}, we obtain the desired result $\mathcal{U}''_{T_0}= \frac{1}{2} \dot\ell^h $.  \endproof

\section{Proof of Theorem \ref{main-super}}

This proof uses essentially the same arguments as the proof of Theorem \ref{mainT},
and for this reason we will skip some details. We suppose that the 
super-Brownian motion $(\mathbf{X}_t)_{t\geq 0}$ is constructed (under the probability measure $\P$) from
a Poisson point measure $\sum_{k\in K} \delta_{\omega_k}$ with intensity
$\alpha\,\N_0$ in the way explained in Section \ref{sec:super}. As previously,
we write $\mathbf{Y}=\int_0^\infty \mathbf{X}_t\,\dd t$ for the total
occupation measure of $\mathbf{X}$. Recall that $(L^x)_{x\in\R}$
is the (continuous) density of $\mathbf{Y}$, and that, for $x\not =0$, $\dot L^x$
stands for the derivative of $y\mapsto L^y$ at $x$, and by
convention $\dot L^0$ is the right derivative at $0$. 

For every $h\geq 0$, we
let $\mathbf{Y}^h_-$, resp. $\mathbf{Y}^h_+$, be the restriction of $\mathbf{Y}$
to $(-\infty,h)$, resp. to $(h,\infty)$, and we write $\wt{\mathbf{Y}}^h_+$
for the pushforward of $\mathbf{Y}^h_+$ under the shift $x\mapsto x-h$. The proof of Theorem \ref{main-super}
then reduces to checking the analog of Proposition \ref{keymainT},
namely the identity
\begin{equation}
\label{key-super}
\E[ F_1(\mathbf{Y}^h_-)\,F_2(\wt{\mathbf{Y}}^h_+)]= \E[F_1(\mathbf{Y}^h_-)\,QF_2(L^h,\frac{1}{2}\dot L^h)],
\end{equation}
where $F_1$ and $F_2$ are nonnegative measurable functions on $\mathcal{M}(\R)$,
and $QF_2$ is defined as in the previous section. 

Consider first the case $h>0$. We note that
$$L^h=\sum_{k\in K} \ell^h(\omega_k)$$
as a consequence of \eqref{occu-super} and the 
fact that there are only finitely many $k\in K$ such that $W^*(\omega_k)\geq h$. 
We can then consider the exit measure
process $(X^h_t)_{t\geq 0}$, which is defined by
$$X^h_t:=\sum_{k\in K} \mathcal{X}^h_t(\omega_k).$$
Note again that there are only finitely many nonzero terms in the right-hand side. 
Then $(X^h_t)_{t\geq 0}$ is (again) a  $\psi$-CSBP, which now starts at
$$X^h_0=Z_h:=\sum_{k\in K} \z_h(\omega_k).$$
We may write $X^h$ as the time change of a L\'evy process $U=(U_t)_{t\geq 0}$ started at $Z_h$, in such a way
that
$$\int_0^\infty X^h_t\,\dd t=T_0:=\inf\{t\geq 0:U_t=0\},$$
and we have 
$$L^h=\sum_{k\in K} \ell^h(\omega_k)=\sum_{k\in K} \int_0^\infty \mathcal{X}^h_t(\omega_k)\,\dd t=\int_0^\infty X^h_t\,\dd t=T_0.$$

There is again a one-to-one correspondence between the jump times
of $X^h$ and the excursions of $\omega_k$ above and below $h$,
for all $k\in K$ (such that $W^*(\omega_k)\geq h$). We can list these excursions
in a sequence $(\omega^{(i)},i\in \N)$ as we did in the preceding section,
and we let $I_-$, resp. $I_+$, be the set of all indices $i$ such that
$\omega^{(i)}$ is an excursion below $h$, resp. above $h$. 
Then, conditionally on the exit measure
process $(X^h_t)_{t\geq 0}$, the excursions $(\omega^{(i)},i\in \N)$ are independent
(and independent of the point measure $\sum_{k\in K} \delta_{\tr(\omega_k)}$), and
the conditional distribution of $\omega^{(i)}$ is $\frac{1}{2}(\N^{*,\delta_i}_h + \check\N^{*,\delta_i}_h)$,
where $\delta_i$ is the jump associated with $\omega^{(i)}$. 

We may then construct the L\'evy processes $U'$ and $U''$ from $U$ in a way exactly similar
as we constructed $\mathcal{U}'$ and $\mathcal{U}''$ from $\mathcal{U}$ in the previous section,
and we have $U'+U''=U$, so that $U'_{T_0}=-U''_{T_0}$.

We can now follow the same route as in the proof of Theorem \ref{mainT}
to arrive at the analog of formula \eqref{key3}, which reads
\begin{align}
\label{key10}
&\E\Bigg[H\Big(\sum_{k\in K} \delta_{\tr_h(\omega_k)}\Big)\,\exp\Big(-\sum_{i\in I_-} A(\omega^{(i)})\Big)\exp\Big(-\sum_{i\in I_+} B(\omega^{(i)})\Big)\Bigg]\\
&\qquad=\E\Bigg[H\Big(\sum_{k\in K} \delta_{\tr_h(\omega_k)}\Big)\,\exp\Big(-\sum_{i\in I_-} A(\omega^{(i)})\Big)\,\Phi_B(T_0,-{U}'_{T_0})\Bigg],\nonumber
\end{align}
with the {\it same} function $\Phi_B$ as in \eqref{key3}. We already know that $T_0=L^h$, and, to complete the proof 
of \eqref{key-super}, we need to verify that $U'_{T_0}=-\frac{1}{2}\,\dot L^h$. This is done by exactly the same method 
as in the proof of Lemma \ref{tech-lemma}, using the approximation
\begin{equation}
\label{conv-deri2}
\frac{1}{\ve^2}\Bigg(\int_h^{h+\ve} L^x\,\dd x - \ve \,L^h\Bigg) \build{\la}_{\ve\to 0}^{} \frac{1}{2} \,\dot L^h
\end{equation}
instead of \eqref{conv-deri}.

Let us consider then the case $h=0$. We start by observing that
$$L^0=\sum_{k\in K} \ell^0(\omega_k).$$
This identity is less immediate than the corresponding one when $h>0$, because the sum
now involves infinitely many nonzero terms, but it is derived in the proof of \cite[Corollary 3.2]{LGR2}. 
We can define the exit measure process $(X^0_t)_{t\geq 0}$ at $0$
by setting $X^0_0=\alpha$ and, for $t>0$,
$$X^0_t=\sum_{k\in K} \mathcal{X}^0_t(\omega_k).$$
As recalled in Section \ref{sec:exitlocal}, we then know that $(X^0_t)_{t\geq 0}$ is a $\psi$-CSBP started at $\alpha$. Moreover,
we have
$$L^0=\sum_{k\in K} \ell^0(\omega_k)=\sum_{k\in K} \int_0^\infty \mathcal{X}^0_t(\omega_k)\,\dd t= \int_0^\infty X^0_t\, \dd t,$$
where the second equality follows from \eqref{local3}. 

As in the case 
$h>0$, there is a one-to-one correspondence between the jump times
of $X^0$ and the excursions of $\omega_k$ above and below $0$,
for all $k\in K$ --- now all $k\in K$ are relevant, but this creates no difficulty, since
$\mathcal{X}^0(\omega_k)$ and $\mathcal{X}^0(\omega_{k'})$ have no common jump
time if $k\not=k'$. 
We can list these excursions
in a sequence $(\omega^{(i)},i\in \N)$ as above.
By a direct application of \cite[Theorem 4]{ALG}, we again get that, conditionally on $(X^0_t)_{t\geq 0}$,
the excursions $(\omega^{(i)},i\in \N)$ are independent, and
the conditional distribution of $\omega^{(i)}$ is $\frac{1}{2}(\N^{*,\delta_i}_0 + \check\N^{*,\delta_i}_0)$,
where $\delta_i$ is the jump of $X^0$ associated with $\omega^{(i)}$. 
The Lamperti time change 
of $X^0$ yields a
L\'evy process $U$ started from $\alpha$, up to time $T_0:=\inf\{t\geq 0:U_t=0\}$, and we can again consider the
L\'evy processes $U'$, resp. $U''$, obtained by ``keeping'' the jumps of $U$ corresponding 
to negative excursions, resp. to positive excursions, and such that $U'_0=\alpha$
and $U''_0=0$. By the same arguments as in the proof of Theorem \ref{mainT}, we
arrive at the analog of \eqref{key10}
(without the term $H(\sum_{k\in K} \delta_{\tr_h(\omega_k)})$ which is now irrelevant). 
Since we already now that $T_0=\int_0^\infty X^0_t\,\dd t=L^0$, it only
remains to verify that $U'_{T_0}=-U''_{T_0}=-\frac{1}{2}\,\dot L^0$. This follows by
a straightforward adaptation of the proof of Lemma \ref{tech-lemma}, using \eqref{conv-deri2} with $h=0$.
This completes the proof of Theorem \ref{main-super}.

\medskip
\rems (i) In the case $h=0$, if instead of using \eqref{conv-deri2}, we consider the approximation
$$\frac{1}{\ve^2}\Bigg(\int_{-\ve}^{0} L^x\,\dd x - \ve \,L^h\Bigg) \build{\la}_{\ve\to 0}^{} -\frac{1}{2} \,\dot L^{0-},$$
the same method leads to the equality $U'_{T_0}-\alpha= -\frac{1}{2} \dot L^{0-}$. Since 
we have also $U'_{T_0}=-U''_{T_0}=-\frac{1}{2}\,\dot L^0$, we get that $\dot L^0=\dot L^{0-}-2\alpha$,
which is consistent with the results of \cite{Sug}. 

\smallskip
\noindent(ii) It is certainly possible to derive \eqref{key-super} more directly from (a stronger form 
of) Proposition \ref{keymainT}. This would still require some technicalities, and
we preferred to use the preceding approach which consists in adapting the proof
of Proposition \ref{keymainT} to a slightly different context. 

\section{Remarks and complements}

\subsection{The transition kernel of $(\ell^x,\dot\ell^x)$}

Our proof of Theorem \ref{mainT} yields a complicated expression for the
transition kernel of the Markov process $(\ell^x,\dot\ell^x)$
(or of the process $(L^x,\dot L^x)$ of Theorem \ref{main-super}). 
First observe that we can use Theorem \ref{re-root-theo}
to verify that $\y(\omega)$ also has a continuously
differentiable density on $(0,\infty)$, $\N^*_0$ a.e. By a scaling argument,
the same holds $\N^{*,z}_0$ a.e. for every $z>0$. 
In other words we can make sense of $(\ell^x,\dot\ell^x)$ for every $x>0$, 
$\N^{*,z}_0$ a.e. For $t>0$ and $y\in\R$, recall the notation 
$U^{\mathrm{br},t,y/2}$
for the bridge of duration $t$ from $0$ to $y/2$ associated with the centered stable L\'evy process with Laplace exponent $\frac{1}{2}\psi$, and
write $(\eta_k)_{k\in\N}$ for the sequence of jumps of $U^{\mathrm{br},t,y/2}$ ranked in nonincreasing order. 

Let $x>0$. Then, under $\N_0$, the law of $(\ell^x,\dot\ell^x)$ knowing that $(\ell^0,\dot\ell^0)=(t,y)$
is the distribution of
$$\Big(\sum_{k\in K} \ell^x(\omega_k),\sum_{k\in K} \dot\ell^x(\omega_k)\Big)$$
where, conditionally on $U^{\mathrm{br},t,y/2}$, the random snake trajectories $\omega_k$
are independent, and, for every $k$, $\omega_k$ is distributed according to $\N^{*,\eta_k}_0$. 
This expression readily follows from Proposition \ref{keymainT}.

We note that there are finitely many nonzero terms in the sums of the last display. To see this, observe
that, for every $z>0$,
$$\N^{*,z}_0(W^*\geq x)=\N^{*,1}_0(W^*\geq x/\sqrt{z})\leq C\,\frac{z^3}{x^6},$$
where $C$ is a constant and the last bound follows from \cite[Corollary 5]{LGR}. Hence,
$$\E\Big[\sum_{k\in K} \mathbf{1}_{\{\ell^x>0\}}\,\Big|\,U^{\mathrm{br},t,y/2}\Big] \leq Cx^{-6}\,\sum_{k\in K}(\eta_k)^3 <\infty,\quad\hbox{a.s.}$$

It would be desirable to obtain a simpler description of the transition kernel of $(\ell^x,\dot\ell^x)$ !

\subsection{Towards a stochastic equation}

The paper \cite{Hon} gives formulas for the local time of a super-Brownian motion 
$(\mathbf{X}_t)_{t\geq 0}$ started at $\delta_0$
and its derivative, in terms of the martingale measure $M$ associated with $(\mathbf{X}_t)_{t\geq 0}$
(see \cite[Section II.5]{Per} for the definition and properties of $M$). 
With our notation, formula (2.11) of \cite{Hon} states that, for every fixed $0<x<y$,
$$\dot L^y - \dot L^x = \int_0^\infty \int (\mathrm{sgn}(x-z)-\mathrm{sgn}(y-z))\,M(\dd z\dd s),$$
where $\mathrm{sgn}(z)=\mathbf{1}_{\{z>0\}}-\mathbf{1}_{\{z<0\}}$. 

Assuming that $(\dot L^x)_{x\geq 0}$ is a semimartingale, one can infer from this representation
that the quadratic variation of its martingale part should be $4\int_0^x L^z\,\dd z$. 
This suggests that $(\dot L^x)_{x\geq 0}$ should satisfy a stochastic differential equation
of the form
$$\dd \dot L^x= 2\sqrt{L^x}\,\dd\beta_x+ h(L^x,\dot L^x)\,\dd x$$
where $(\beta_x)_{x\geq 0}$ denotes a linear Brownian motion, and $h$ 
is a measurable function 
on $\R_+\times \R$. This is indeed the case as the
forthcoming work \cite{LGP} will show. 

The equation in the last display is very
close to the one that is conjectured to hold for the density of ISE in \cite{CM}.
Note however that the (conjectured) drift term in \cite[Conjecture 7]{CM} is 
more complicated because of the conditioning involved in the definition of ISE.

As a final remark, the stochastic equation in the last display is of course 
reminiscent of the equation $\dd X_x=2\sqrt{X_x}\,\dd\beta_x$ which (by
the Ray-Knight theorems) holds if $X_x$ is the local time at level $x>0$ 
of a positive Brownian excursion distributed according to the It\^o measure.

\subsection{Brownian geometry}

The Brownian sphere, or Brownian map, is a random measure metric space
$(\bm,D,\mathrm{vol})$ that arises
as the scaling limit in the Gromov-Hausdorff sense of many different classes 
of random planar maps (see in particular \cite{Uniqueness,Mie}). The Brownian 
sphere is constructed as the quotient space $\bm=\t/\approx$ of the Brownian 
tree $\t$ for an equivalence relation $\approx$ defined in terms of the labels $(V_a)_{a\in\t}$,
and the volume measure on $\bm$ is just the pushforward of the volume
measure $\mathrm{Vol}$ on $\t$ under the canonical projection. Under $\N_0(\cdot\mid \sigma=1)$, we speak of the standard Brownian sphere
(with total volume equal to $1$),
but it is also of interest to consider the ``free'' Brownian sphere defined under $\N_0$.
The equivalence relation $\approx$ is such that we have $V_a=V_{a'}$ whenever
$a$ and $a'$ are two points of $\t$
such that $a\approx a'$. Thanks of this property, one can make
sense of the label $V_\bx$ for any point $\bx$ of $\bm=\t/\approx$. 

The Brownian sphere comes with two distinguished points, namely $\bx_0$, which is the
equivalence class of the root of $\t$, and $\bx_*$,
which is the
equivalence class of the point of $\t$ with minimal label
(in a sense that can be made precise, these two points are uniformly
distributed over $\bm$). Moreover, we have $D(\bx_*,\bx)= V_\bx-V_{\bx_*}$ for
every $\bx\in \bm$: up to a shift, labels correspond to distances from the
distinguished point $\bx_*$. The next proposition is then a straighforward consequence
of the preceding results. To simplify notation, we write $m_*=-V_{\bx_*}=D(\bx_0,\bx_*)$.

\begin{proposition}
\label{geo-brown}
For every $r\geq 0$, let $\v_r$ be the volume of the closed ball of radius $r$ 
centered at $\bx_*$ in the Brownian sphere $\bm$. Then, $\N_0$ a.e. the
function $r\mapsto \v_r$ is twice continuously differentiable on $[0,\infty)$, and
we denote its first and second derivative by $\v'_r$ and $\v''_r$.
Moreover, the random process $(\v_{m_*+r},\v'_{m_*+r},\v''_{m_*+r})_{r\geq 0}$
is time-homogeneous Markov under $\N_0$.
\end{proposition}

\proof
By the definition of the volume measure on $\bm$, and the formula for distances from $\bx_*$,
$$\v_r=\mathrm{Vol}(\{a\in\t:V_a\leq r-m_*\})=\int_{-\infty}^{r-m_*} \,\ell^x\,\dd x.$$
From the fact that $x\mapsto\ell^x$ is continuously differentiable, we 
thus get that the mapping $r\mapsto \v_r$ is twice continuously differentiable,
and moreover $\v'_{m_*+r}=\ell^r$ and $\v''_{m_*+r}=\dot\ell^r$. Then we just
have to apply Theorem \ref{mainT}. \endproof

\smallskip
Informally, $\v'_r$ represents the ``area'' of the sphere $\{x\in\bm:D(\bx_*,\bx)=r\}$.
Furthermore, Lemma \ref{ident-derivative} allows us to interpret
$\v''_r$ as twice the (renormalized) sum of the boundary sizes of connected components
of the complement of the closed ball of radius $r$ centered at $\bx_*$: in the canonical projection 
from $\t$ onto $\bm$, these connected components correspond to
the excursions above level $r-m_*$ (see the beginning of \cite[Section 12]{Disks}).

%
%

\begin{acks}[Acknowledgments]
This work was motivated by a very stimulating lecture of Guillaume Chapuy
at the CIRM Conference on Random Geometry in January 2022. I thank Guillaume Chapuy and 
Jean-Fran\c cois Marckert for keeping me informed of their on-going work \cite{CM}. I also thank
Lo\"ic Chaumont for letting me know about reference \cite{BCP}. It is
a pleasure to thank Jieliang Hong and Ed Perkins for a useful conversation
at the 2022 PIMS Summer School in Probability. Finally, I thank two anonymous referees 
for their careful reading of the manuscript and, in particular, for pointing at an error 
in the initial proof of Lemma \ref{tech-lemma}.
\end{acks}
\begin{funding}
This work was supported by the ERC Advanced Grant 740943 {\sc GeoBrown}.
\end{funding}


\begin{thebibliography}{99}


\bibitem{ALG}
{C. Abraham, J.-F. Le Gall}, Excursion theory for Brownian motion
indexed by the Brownian tree. {\it J. Eur. Math. Soc. (JEMS)} 20, 2951--3016 (2018)



\bibitem{Ald}
{D. Aldous}, The continuum random tree I. {\it Ann. Probab.}, 19, 1--28 (1991)


\bibitem{Al2}
{D. Aldous}, Tree-based models for random distribution of mass. {\it J. Statist. Phys.} 73, 625--641  (1993)

\bibitem{Be}
{J. Bertoin}, L\'evy processes. Cambridge University Press, 1996.

\bibitem{BCP}
{J. Bertoin, L. Chaumont, J. Pitman}, Path transformations of first passage bridges. {\it Electron. Comm. Probab.} 8, 155--166 (2003)

\bibitem{BM0}
{M. Bousquet-M\'elou}, 
Limit laws for embedded trees: applications to the integrated superBrownian excursion. {\it Random Structures Algorithms} 29, 475--523 (2006)


\bibitem{BMJ}
{M. Bousquet-M\'elou, S. Janson}, The density of the ISE and local limit laws for 
embedded trees. {\it Ann. Appl. Probab.} 16, 1597--1632 (2006)


\bibitem{BCL}
M. Bramson, J.T. Cox, J.-F. Le Gall, Super-Brownian limits of voter model clusters. {\it Ann. Probab.} 29, 1001--1032 (2001)


\bibitem{CM}
{G. Chapuy, J.-F. Marckert}, Note on the density of ISE and a related diffusion. Preprint, arXiv:2210.10159


\bibitem{CLG} N. Curien, J.-F. Le Gall, The hull process of the Brownian plane.
{\it Probab. Theory Related Fields} 166, 187--231 (2016)

\bibitem{DS}
{E. Derbez, G. Slade}, The scaling limit of lattice trees in high dimensions. 
{\it Comm. Math. Phys.} 193, 69--104 (1998)
%
%
%
%
%


\bibitem{FPY}
{P. Fitzsimmons, J. Pitman, M. Yor}, Markovian bridges: construction, Palm interpretation and splicing. Seminar on Stochastic Processes 1992, 101--134, Progr. Probab., 33, Birkh\"auser Boston, 1993. 

\bibitem{HaS}
{T. Hara, G. Slade}, The scaling limit of the incipient infinite cluster in high-dimensional percolation. II. 
Integrated super-Brownian excursion. Probabilistic techniques in equilibrium and nonequilibrium statistical physics. 
{\it J. Math. Phys.} 41, 1244--1293 (2000)


\bibitem{Hon}
{J. Hong}, Improved H\"older continuity near the boundary of one-dimensional 
super-Brownian motion. {\it Electron. Comm. Probability} 24, no 28, 1--12 (2019)

\bibitem{King}
J.F.C. Kingman, Poisson processes. Oxford University Press, 1993.


\bibitem{Lam}
J. Lamperti, 
Continuous state branching processes. 
{\it Bull. Amer. Math. Soc.} 73, 382--386 (1967)


\bibitem{Zurich} {J.-F. Le Gall}, {\it Spatial branching processes, random snakes and 
partial differential equations}. {Lectures in Mathematics ETH Z\"urich}. Birkh\"auser, Boston, 1999.

 \bibitem{probasur} 
 {J.-F. Le Gall}, {Random trees and applications}.
{\it Probab. Surveys} {2}, 245--311 (2005)


\bibitem{Uniqueness} {J.-F. Le Gall}, Uniqueness and universality of the Brownian map.
{\it Ann. Probab.} {41}, 2880--2960 (2013)

\bibitem{subor} 
{J.-F. Le Gall}, Subordination of trees and the Brownian map.
{\it Probab. Theory Related Fields} {171}, 819--864 (2018)

\bibitem{Disks} 
J.-F. Le Gall, Brownian disks and the Brownian snake. {\it Ann. Inst. H. Poincar\'e Probab. Stat.} 55, 237--313 (2019)

\bibitem{LGP}
J.-F. Le Gall, E. Perkins, A stochastic differential equation for local times
of tree-indexed Brownian motion. In preparation. 

%
%


\bibitem{LGR}
{J.-F. Le Gall, A. Riera}, Growth-fragmentation processes in Brownian motion indexed by the Brownian tree. {\it Ann. Probab.} 48, 1742--1784 (2020)

\bibitem{LGR2}
{J.-F. Le Gall, A. Riera}, Some explicit distributions for Brownian motion indexed by the Brownian tree. {\it Markov Processes Relat. Fields} 26, 659--686 (2020)

\bibitem{Spine}
{J.-F. Le Gall, A. Riera}, Spine representations for non-compact models of random geometry. {\it Probab. Th. Rel. Fields} 181, 571-645 (2021)

\bibitem{LGW}
{J.-F. Le Gall, M. Weill}, Conditioned Brownian trees. {\it Ann. Inst. Henri Poincar\'e Probab. Stat.} 42, 455--489 (2006)

\bibitem{Mie}
{G. Miermont},
The Brownian map is the scaling limit of uniform random plane quadrangulations.
{\it Acta Math.}, {210}, 319--401 (2013)


\bibitem{Per}
{E.A. Perkins}, Dawson-Watanabe superprocesses and measure-valued diffusions. 
Ecole d'\'et\'e de probabilit\'es de Saint-Flour 1999. {\it Lecture Notes Math.} 1781. Springer 2002

\bibitem{Pit}
{J. Pitman}, Combinatorial stochastic processes.
Ecole d'\'et\'e de probabilit\'es de Saint-Flour 2002. {\it Lecture Notes Math.} 1875. Springer 2006

\bibitem{Sug}
{S. Sugitani},
Some properties for the measure-valued branching diffusion processes. 
{\it J. Math. Soc. Japan} 41, 437--462 (1989)

\bibitem{Yor}
{M. Yor.} On some exponential functionals of Brownian motion. {\it Adv. Appl. Prob.} 24, 509--531 (1992)


\end{thebibliography}
\end{document}